 \documentclass[11pt]{article}

 \usepackage{amsfonts}
 \usepackage{amssymb}
 \usepackage{amsxtra}

 \newcommand{\obra}[3]{{\sc #1} {\em #2}. {#3}.}

 \newtheorem{theorem}{\bf Theorem}
 \newtheorem{lemma}[theorem]{\bf Lemma}
 \newtheorem{proposition}[theorem]{\bf Proposition}
 \newtheorem{definition}[theorem]{\bf Definition}
 \newtheorem{corollary}[theorem]{\bf Corollary}
 \newtheorem{remark}[theorem]{Remark}
 
 \newtheorem{example}[theorem]{Example}

\newcommand{\deuxbis}{2.4'}

  \newenvironment{proof}{{\em Proof \/.-}}
   {\hfill $\square$\newline}

 \newcommand{\R}{\mathbb{R}}
 
 \newcommand{\N}{\mathbb{N}}
 \newcommand{\Z}{\mathbb{Z}}
  \newcommand{\C}{\mathbb{C}}

 \newcommand{\FF}{\mathcal{F}}
 \newcommand{\CC}{\mathcal{C}}
  
   \newcommand{\RR}{\mathcal{R}}
   \newcommand{\DD}{\mathcal{D}}
 \newcommand{\AAA}{\mathcal{A}}
 \newcommand{\SC}{\mathcal{S}}
\newcommand{\matr}[2]{{\left(\begin{array}{#1} #2 \end{array}\right)}}
 \newcommand{\xunoxm}{(x_1,\ldots,x_m)}
 \newcommand{\xx}{{\bf x}}
 \newcommand{\yy}{{\bf y}}
 \newcommand{\zz}{{\bf z}}
  \newcommand{\ww}{{\bf w}}
\newcommand{\ra}{\rightarrow}
\newcommand{\eps}{\varepsilon}
\newcommand{\val}{\mbox{\rm val\,}}
\renewcommand{\Re}{\mbox{\rm Re}\,}
\newcommand{\labe}[1]{\label{#1}}

\title{Quasi-analytic solutions of analytic ordinary differential equations
 and o-minimal structures}
\author{{\sc Rolin, J.-P.\footnote{Laboratoire de Topologie UMR-5584, Universit\'{e} de Bourgogne, BP138, 21004 Dijon cedex (France). Fax: (33) 3 80 39 58 99. E-mail: rolin@u-bourgogne.fr}; Sanz, F.\footnote{Departamento de \'{A}lgebra, Geometr\'{\i}a y Topolog\'{\i}a, Universidad de Valladolid, Prado de la Magdalena, s/n, E-47005, Valladolid (Spain). Fax: (34) 983 42 37 88. E-mail: fsanz@agt.uva.es}; Sch\"{a}fke, R.\footnote{Institut de Recherche Math\'{e}matique Avanc\'{e}e (IRMA), Universit\'{e} Louis Pasteur de Strasbourg, 7, Ren\'{e} Descartes
67084, Strasbourg cedex (France). Fax: (33) 3 90 24 03 28. E-mail: schaefke@math.u-strasbg.fr}}}
\date{\today}
 
\begin{document}
 \maketitle
\begin{abstract}
It is well known that the non-spiraling leaves of real
analytic foliations of codimension 1 all belong to the
same o-minimal structure.
Naturally, the question arises if the same statement is true for
non-oscillating trajectories of real analytic vector fields.

We show, under certain assumptions, that such a trajectory
generates an o-minimal and model complete structure together with the
analytic functions.
The proof uses the asymptotic theory of irregular singular ordinary
differential equations in order to establish a quasi-analyticity result
from which the main theorem follows.

As applications, we present an infinite family of o-minimal
structures such that any two of them do not admit a common extension,
and we construct a non-oscillating trajectory of a real analytic vector
field in $\R^5$ that is not definable in any o-minimal extension of $\R$.
\end{abstract}

\section{Introduction}
\setcounter{theorem}{0}\setcounter{equation}{0}

 Consider a
 system of ordinary differential equations having an irregular singular point:
 \begin{equation}\labe{sistema}
 x^{p+1}\frac{d\yy}{dx}=A(x,\yy), \;\; \yy\in\R^r,
 \end{equation}
 where $A$ is real analytic in a neighborhood of $0\in\R^{r+1}$
 with $A(0,0)=0$ and the {\em Poincar\'{e} rank} $p$ is greater or equal to one.
In a series of articles
\cite{Lio-Mil-Spe},
\cite{Can-Mou-Rol},
\cite{Can-Mou-San}, {\em non-oscillating solutions}
of (\ref{sistema}) have been
studied; here a solution (vector) $H=(H_1,...,H_r):(0,\eps]\rightarrow\R^r$
of (\ref{sistema}) tending to 0  as $x\rightarrow 0$ is called non-oscillating
if for any function $f$ real analytic in some neighborhood of
$0\in\R^{r+1}$, the function $x\mapsto f(x,H(x))$ either is identically 0
or has finitely many zeros.

One of the motivations of our work is the more general question whether
the {\em structure} $\R_{an,H}$
generated by $H_1,...,H_r$ and the restricted analytic functions is {\em o-minimal},
i.e.\ the sets definable in that structure have a finite number of
connected components (with respect to the usual topology of $\R^m$).
The definable sets of $\R_{an,H}$ are those that belong
to the
smallest family of subsets of $\R^m$,
$m\in\N$, that contains the graphs of addition, multiplication,
restrictions of analytic functions to cartesian products of compact intervals and $H_1,...,H_r$
and that is
closed under finite unions and intersections, cartesian products, complements
and linear projections.
It is easy to show that the graphs of the above functions $f(x,H(x))$,
$f$ real analytic in a
neighborhood of 0, are definable, but more generally also the graph
of any function
obtained from restricted analytic functions and
$H_j, \,j=1,...,r$ by (repeated) compositions and resolution of
implicit equations. The desired o-minimality of  $\R_{an,H}$ thus would
imply that they all have a finite number of zeros.

For bounded non-spiraling leaves of real analytic foliations of codimension 1,
such a result has been shown (see \cite{Wil}, \cite{Lio-Rol}): the structure generated by all such leaves is o-minimal.

In our main  result (Theorem \ref{main3}), we consider a system (\ref{sistema}) such that
a) the eigenvalues
$\lambda_1,\ldots,\lambda_r$ of $\frac{\partial
A}{\partial\yy}(0,0)$ are nonzero and have distinct arguments;
and b) the Stokes phenomenon of the formal solution $\widehat H$ is non-trivial
for at least one direction corresponding to each $\lambda_j,\,j=1,...,r$.
Observe that this implies the divergence of $\widehat H$.
We consider an arbitrary solution
$H:(0,\varepsilon]\rightarrow\R^r$ of such an equation having an
asymptotic expansion $H(x)\sim\widehat{H}(x)=\sum_{n=1}^\infty h_n\,x^n$
as $x\rightarrow0$ and prove that the above structure $\R_{an,H}$
generated by the analytic functions and $H_1,...,H_r$ is
o-minimal and model complete;
a structure is called {\em model complete}, if the same family of subsets is
generated without taking complements.
Observe that our theorem can also be applied to $C^\infty$-solutions $H:[-\varepsilon,\varepsilon]\rightarrow\R^r$ of (\ref{sistema})
if $p$ is odd:
it suffices to consider $\tilde H:(0,\varepsilon]\rightarrow\R^{2r}$ defined
by $\tilde H(x)= (H(x),H(-x))$. Condition a) has to be replaced by the following one: the eigenvalues
$\lambda_1,\ldots,\lambda_r$ of $\frac{\partial
A}{\partial\yy}(0,0)$ are nonzero and have arguments distinct mod $\pi\Z$.

The model completeness is linked to the problem of describing $\R_{an,H}$
in terms of $H$ and analytic functions: here we show (Proposition \ref{H-conjuntos})
that any definable set
in $\R_{an,H}$ can be obtained as finite union of linear projections
of sets of the form
$A=X\cap V_1\cap V_2\cdots\cap V_n$
where $X$ is globally semi-analytic (see subsection 2.1)
 and $V_i=\{\xunoxm\;/\;x_k=H_s(x_l)\}$ for some $1\leq k,l\leq m$ and
 $1\leq s\leq r$.

As applications of our result, we present a two-parameter family of solutions of
some two-dimensional system (\ref{sistema}) such that the o-minimal structures
generated by each of them are mutually incompatible, i.e.\ any two of them do
not admit a common o-minimal extension -- this is very different from the theory of
(bounded) non-spiraling leaves of analytic foliations: they {\em all} belong to one common o-minimal structure. Moreover, we construct a non-oscillating
trajectory of a system (\ref{sistema}) in $\R^4$ that cannot belong to any
o-minimal structure.

Our main result is a consequence of two theorems. The first
(Theorem~\ref{main1}) states that strong quasi-analyticity of $H$
(see Definition~\ref{def:QA}) implies o-minimality and model completeness of
the structure $\R_{an,H}$, the second (Theorem~\deuxbis{}) states that our
conditions a), b) imply strong analytic transcendence of the formal solution
$\widehat H(x)$ of (\ref{sistema}) which in turn yields that the above $H$
is strongly quasi-analytic. Its proof is based on the theory of irregular singular points of ordinary differential equations in the complex domain.

The article is organized as follows: in section 2, we present precise
definitions and theorems as well as examples and applications, section 3 is
devoted to the proof of Theorem~\ref{main1}, section 4 to that of Theorem~\deuxbis{}. An outline of these proofs can be found at the end of
subsection 2.2.

\section{Results}
\setcounter{theorem}{0}\setcounter{equation}{0}

\subsection{Preliminaries}
Before stating the results, let us briefly recall the terminology of o-minimal
structures (see
 \cite{vdD} for an introduction).
Consider a family $\mathcal{F}=(\FF_m)_{m\in\N}$ of sets $\mathcal{F}_m$ of functions
$f:\R^m\rightarrow\R$. A subset of
 $\R^m$ is said to be {\em definable in the structure $\R_\FF$ generated
 by $\FF$} (or {\em in the extension of the real field by $\FF$})
 if it belongs to the smallest collection of subsets of $\R^m,\,m\in\N,$ satisfying:
 \begin{itemize}
 \item[S1)] It is closed under cartesian products,
  finite unions and intersections, complements and under images of linear projections
  $\R^{m+1}\rightarrow\R^m$.
 \item[S2)] It contains the diagonals
$\Delta_{ij}=\{\xunoxm\in\R^m\;/\;x_i=x_j\}$ for $1\leq i<j\leq
m$.
 \item[S3)] It contains the point sets $\{a\}$ for $a\in\R$,
 as well as the graphs of addition
 and multiplication $+,\cdot:\R^2\rightarrow\R$.
 \item[S4)] It contains the graph of any function in $\mathcal{F}$.
 \end{itemize}

If some function of the family $\cal F$ is defined on a proper subset $D\subset \R^m$,
then it is understood without mentioning this, that the function is extended to
$\R^m$ by the value 0 in the above definition.

 The structure $\R_{\mathcal{F}}$ is called {\em o-minimal} if
 every definable set has a finite number of
 connected components (with respect to the ordinary topology).
It is called {\em model complete} if in the above
definition the condition of closure under complements is superfluous.

For $\FF=\emptyset$, the structure $\R_{\FF}=\R_{alg}$ is o-minimal
and model complete and its definable sets are the {\em
semialgebraic
 sets}, i.e.\ finite unions of subsets of $\R^m$ of the form
 $\{\xx\in\R^m\;/\;P(\xx)=0,Q_1(\xx)>0,\ldots,Q_s(\xx)>0\}$ where $P,Q_1,\ldots,Q_s$
 are polynomials. Another
 classical example uses real
 analytic functions instead of polynomials.
 A subset $A\subset\R^n$ is called {\em
 semianalytic at a point $a\in\R^n$} if there exists an open
 neighborhood $U$ of $a$ such that $A\cap U$ is a finite union of
 sets of the form $\{\xx\in
 U\;/\;f(\xx)=0,g_1(\xx)>0,\ldots,g_s(\xx)>0\}$, where
 $f,g_1,\ldots,g_s$ are analytic functions on $U$. The set
 $A\subset\R^n$ is called {\em global semianalytic} if $\rho_n(A)\subset [-1,1]^n$
 is semianalytic at every point, where
 $\rho_n(x_1,\ldots,x_n)=(x_1/\sqrt{1+x^{2}_{1}},\ldots,x_n/\sqrt{1+x^{2}_{n}})$.
Finally,
  a
 {\em global subanalytic set} of $\R^m$ is the linear projection of
 a global semianalytic set $A\subset\R^n$ for some $n\geq m$.
 The collection of global subanalytic sets
 is precisely the collection of definable sets in the
 structure $\R_{an}$ generated by the set $\widetilde{\FF}_{an}$ of {\em
 restricted analytic functions}.
These are the functions $\tilde{f}:\R^m\rightarrow\R$ for which
 there exists
 an analytic function $f$ defined in a neighborhood of $[-1,1]^m$ such that
 $\tilde{f}=f$ on $[-1,1]^m$
 and $\tilde{f}=0$ outside $[-1,1]^m$. It is known that $\R_{an}$
 is o-minimal and model complete: this result is
 a consequence of Gabrielov's
 Complement Theorem asserting that the complement of a global subanalytic set
 is again global subanalytic \cite{Gab}.

\subsection{Statement of the results}
Consider a system (\ref{sistema}) of analytic  ordinary differential equations and
a solution (vector) $H=(H_1,\ldots,H_r):(0,\varepsilon]\rightarrow\R^r$ with $\lim_{x\rightarrow 0}H(x)=0$
admitting an asymptotic expansion
$H(x)\sim\widehat{H}(x)=\sum_{n=0}^\infty h_n x^n\in\R[[x]]^r$
as $\R_+\ni x\ra 0$, i.e.\ for any natural number $N\geq 0$,
there are constants $C_N,\delta_N>0$ such that
$$
\parallel H(x)-\sum_{n=0}^N h_n x^n\parallel\leq C_Nx^{N+1},\;0<x<\delta_N.
$$
We prove below that $\widehat{H}(x)$ is a formal
solution of (\ref{sistema}) and that $H$ extends to a $\CC^\infty$
function on $[0,\eps]$ having $\widehat{H}(x)$ as its Taylor
series at the origin (cf.\ subsection 3.1 for details).

As a matter of notation, if $\phi(x)$ is a formal power series or
a $\CC^\infty$ function at the origin,
let $J_k\phi\,(x)=\sum_{i=0}^{k}\frac{\phi^{(i)}(0)}{i!}\,x^i\in\R[x]$
denote the $k$-jet of $\phi$ as $\R^+\ni x\ra0$ and let $T_k\phi\,(x)$ denote
$(\phi(x)-J_k\phi\,(x))/x^{k}$. Denote by $\val\phi$ the valuation of $\phi$
at the origin, i.e.\ the minimum of all $l$ with $\phi^{(l)}(0)\neq 0$ and by
$\deg P$ the degree of a polynomial $P$.
\begin{definition}\labe{def:QA}
We say that a solution  $H=(H_1,...,H_r)$ of (\ref{sistema}) is {\em strongly quasi-analytic}
if it tends to 0 at the origin, admits an asymptotic expansion
$H(x)\sim \widehat H(x)$ as $x\rightarrow 0$
and satisfies the following condition:

\vspace{.2cm}  (SQA) \;\;\; If $k\geq 0$, $n\geq 1$, an
analytic function $f\in\R\{x,z_{11},\ldots,z_{rn}\}$ with $f(0)=0$
and polynomials $P_1(x),\ldots,P_n(x)$ with $\val P_l >0$ and $P_l^{(\val P_l)}(0)>0$
are given, then one has
$$
f(x,\{T_k\widehat{H}_j\,(P_l(x))\}_{j,l})\equiv 0\; \Longrightarrow\; f(x,\{T_k
H_j\,(P_l(x))\}_{j,l})\equiv 0.
$$
\end{definition}

Here, the word ``strongly" reflects the fact that
polynomials $P_l$ are allowed as arguments of $H$;
this will be necessary in the sequel (see also Example
\ref{ex:euler2} below). The condition
on the derivatives of the $P_l$ is imposed
by the fact that $H$ is defined only for $x>0$.

The condition (SQA) not only excludes exponentially small solutions of (1)
(if those exist as e.g.\ for $x^2y'=y$), but also the possibility that exponentially small
functions might be obtained as certain combinations of $H_j$ and analytic functions.

\begin{theorem}\labe{main1}
Suppose that $H:(0,\varepsilon]\rightarrow\R^r$ is a strongly
quasi-analytic solution of the system (\ref{sistema}).
Then the structure
$\R_{an,H}:=\R_{\widetilde{\FF}_{an}\cup\{H\}}$
generated by the restricted analytic functions and ${H}$
is o-minimal and model complete.
 \end{theorem}
We can be more explicit about the description of the
definable sets in
 $\R_{an,H}$.
Let us say that a set $A\subset\R^m$ is {\em $H$-semianalytic at a
point $a\in\R^m$} if there exists $\delta>0$ such that $(A-a)\cap
[-\delta,\delta]^{m}$ is a finite union of sets of the form
$X\cap V_1\cap\cdots\cap V_n$, for $n\geq 0$,
where $X$ is semianalytic at every point of $[-\delta,\delta]^{m}$
and $V_i=\{\xunoxm\in
[-\delta,\delta]^{m}\,\mid \,0<x_l\leq\eps,\,x_k=H_s(x_l)\}$ for some
$1\leq k,l\leq m$ and $1\leq s\leq r$. The set $A$ is said to be
{\em global $H$-semianalytic} if $\rho_m(A)$ is $H$-semianalytic
at any point of $[-1,1]^m$, where
 $\rho_m(x_1,\ldots,x_m)=(x_1/\sqrt{1+x^{2}_{1}},\ldots,x_m/\sqrt{1+x^{2}_{m}})$.
 Finally, a set
$B\subset\R^m$ is said to be {\em (global) $H$-subanalytic} if
there exists $n\geq m$ and some global $H$-semianalytic set
$A\subset\R^n$ such that $B=\pi(A)$, where
$\pi:\R^n\rightarrow\R^m$ is the projection onto the first $m$
coordinates.

\begin{proposition}\labe{H-conjuntos}
Under the hypotheses of Theorem \ref{main1}, the
collection of $H$-sub\-ana\-ly\-tic sets is stable with respect to
taking complements
and coincides with the collection of definable sets in
$\R_{an,H}$.
\end{proposition}

Property (SQA) is
trivially satisfied for its sum if the formal series solution
$\widehat{H}(x)$ converges. In this case, the sum of
$\widehat{H}(x)$ is the only strongly quasi-analytic solution
among those having $\widehat{H}(x)$ as an asymptotic expansion
(irrespective of the fact that they may generate an o-minimal
structure: $H(x)=\exp(-1/x)$, $x>0$, or not:\\
$H(x)=(H_1(x),H_2(x))=(\exp(-1/x)\sin(1/x),\exp(-1/x)\cos(1/x))$,
$x>0$).

If $\widehat{H}(x)$ is {\em multisummable} in the positive
real direction $\R^+$ in the sense of \cite{Ram}, then the restriction
$H:(0,\varepsilon)\rightarrow\R$ of the multisum of $\widehat{H}$
along $\R^+$ is strongly quasi-analytic. (We point out
that, for this particular solution, the o-minimality of the
structure $\R_{an,H}$ is a consequence of the results in
\cite{vdD-Spe}).

In the above examples only particular solutions are
(or are shown to be) strongly quasi-analytic. We will state conditions on the system
(\ref{sistema}) which imply (SQA) for {\em any} solution
$H:(0,\varepsilon]\rightarrow\R^r$ having an asymptotic expansion
$H(x)\sim\widehat{H}(x)$ as $\R^+\ni x\ra0$:

$\bullet$ We assume that the linear part $A_0=\frac{\partial
A}{\partial \yy}(0)$ is non singular and that its eigenvalues
$\{\lambda_1,\ldots,\lambda_r\}$ satisfy:
\begin{equation}\labe{eq:argumento-autovalores}
{\rm arg}(\lambda_i)\not\equiv{\rm arg}(\lambda_j)\;\;{\rm
mod}\;2\pi\Z\;\;\;\;\mbox{ if }i\neq j.
\end{equation}

$\bullet$ Condition (\ref{eq:argumento-autovalores}) implies that the
system has a unique formal power series solution
$\widehat{H}(x)\in\R[[x]]^r$ (see for instance \cite{Was}).
Moreover, $\widehat{H}(x)$ is {\em $p$-summable} (see subsection~4.1 for details) in each direction
$d_\theta=\{z\;/\;{\rm arg}(z)=\theta\}$ in the complex plane except
for the {\em singular directions}
$$\{d_{\theta_{l,j}}\,/\,p\,\theta_{l,j}={\rm arg}(\lambda_j)+2\pi l\}_
{1\leq j\leq r;\,0\leq l\leq p-1}.
$$
 To each
$d_{\theta_{l,j}}$ corresponds a {\em Stokes coefficient}
$c_{l,j}\in\C$
(see again subsection 4.1 for details).
If $\widehat{H}(x)$ diverges, at least one of the Stokes coefficients must be nonzero; we assume that moreover
\begin{equation}\labe{stokes}
\forall j\in\{1,\ldots,r\}\;\;\exists\,l=l(j)\in\{0,\ldots,p-1\}\mbox{
with }c_{l,j}\neq 0.
\end{equation}
Observe that $\widehat H(x)$ is not summable in the real positive direction
and thus the results of \cite{vdD-Spe} do not apply
if $c_{l,j}\neq0$ for some $(l,j)$ with $\theta_{l,j}\equiv0\;\;{\rm mod\;2\pi\Z}$.
\begin{theorem}\labe{main2}
 Consider a system of analytic ordinary differential equations (\ref{sistema}) such that conditions
 (\ref{eq:argumento-autovalores}) and (\ref{stokes}) are fulfilled.
 Then any  solution
$H=(H_1,\ldots,H_r):(0,\varepsilon]\rightarrow\R^r$ of
 (\ref{sistema}) having an asymptotic expansion $H(x)\sim\widehat{H}(x)$ as $\R^+\ni x\ra0$
 is strongly quasi-analytic.
 \end{theorem}
In fact, we will prove the following stronger result: 

\vspace{.2cm}
 \noindent {\bf Theorem \deuxbis{}} {\em Under the hypothesis of
Theorem \ref{main2}, the formal solution $\widehat{H}$ is {\em strongly analytically transcendental} in the following sense:

\vspace{.3cm} \!\!\! (SAT) \;\;\; If $k\geq 0$, $n\geq 0$, an
analytic function $f\in\R\{x,z_{11},\ldots,z_{rn}\}$ with $f(0)=0$
and distinct polynomials $P_1(x),\ldots,P_n(x)$
with $\deg P_l < (p+1) \val P_l $ and $P_l^{(\val P_l)}(0)>0$ are given, then one has
$$
f(x,\{T_k\widehat{H}_j\,(P_l(x))\}_{j,l})\equiv 0\; \Longrightarrow\; f\equiv
0.
$$}

\vspace{.2cm} \noindent We will prove below (cf.\ Lemma
\ref{lema:SATimplicaSQA}) that property (SAT) for the
formal solution $\widehat{H}(x)$ of (\ref{sistema}) implies property (SQA) for any
solution $H(x)$ with $H(x)\sim\widehat{H}(x)$ as $\R^+\ni x\ra0$. Combining
Theorems \ref{main1} and \ref{main2}, we can state our main result
\begin{theorem}\labe{main3}
 Let $x^{p+1}d\yy/dx=A(x,\yy)$ be a system of analytic ordinary differential equations
satisfying the above conditions (\ref{eq:argumento-autovalores}) and (\ref{stokes}).
 Then,
 for any solution $H:(0,\varepsilon]\rightarrow\R^r$
admitting an asymptotic expansion $H(x)\sim\widehat{H}(x)$ as $\R^+\ni x\ra0$, the
 structure
$\R_{an,H}$ generated by the restricted analytic functions and $H$
is o-minimal and model complete.
 \end{theorem}

\par
The proof of these results is organized as follows. In paragraph
3, we prove Theorem \ref{main1} and Proposition \ref{H-conjuntos};
the main idea is to apply a fundamental result of \cite{Rol-Spe-Wil}
to a certain class of functions corresponding to $H$.
Let us give here a brief description of the principal arguments.

First, we show that we can suppose the solution $H$ to be a
$\CC^\infty$ function on $[-\varepsilon,\varepsilon]$, by
considering the ramification $x\mapsto H(x^2)$. Then we consider
the smallest class $\AAA=(\AAA_m)_{m\geq 1}$ of germs (at the origin of $\R^m$) of
$\CC^\infty$ functions  that contains the
germs of analytic functions and those of the components $H_j$ of
$H$ and that is stable under composition, division by
monomials and under taking solutions of implicit
equations.
It is shown in \cite{Rol-Spe-Wil}, that the quasi-analyticity
of the class $\AAA$ suffices to prove a
{\em Theorem of the Complement for $\AAA$}: the
``subanalytic sets" defined by using elements of $\AAA$ have a
finite number of connected components and they form a family
that is stable under taking complements. Theorem \ref{main1} and Proposition
\ref{H-conjuntos} follow by classical arguments.

The difficult part is to show that $\AAA$ is a {\em
quasi-analytic} class: if $\phi\in\AAA_m$ has a vanishing Taylor
series at $0\in\R^m$ then $\phi=0$. Here, we first reduce
the problem, by restricting to appropriate paths, to the
quasi-analyticity of the subclass $\AAA_1$ of germs of functions of a single
variable.

The quasi-analyticity of $\AAA_1$ is not an obvious consequence
of hypothesis
(SQA); $\AAA_1$ seems to be a very large class of functions,
containing for instance compositions $H_j\circ H_l$
of components of $H$ and implicitly defined functions.
By means of the crucial
Lemma \ref{morfismo} and a description of the functions of $\AAA$
(Lemma \ref{lema2}),
 we get rid of this ``double
transcendence" and we can prove that $\AAA_1$
is the set of so-called {\em simple functions}: germs of
functions of the form $x\mapsto f(x,\{T_k H_j\,(P_l(x))\}_{j,l})$
where $f$ is analytic and $P_l$ are polynomials. This reduction
justifies our condition (SQA).

In paragraph 4, we prove Theorem \deuxbis{}; roughly, we proceed as follows:
since $T_k H(x)$ satisfies a system
like (\ref{sistema}) as well,
we can first suppose that $k=0$. Consider the formal series
$\widehat{F}(x)=f(x,\{\widehat{H}_j(P_l(x))\})$,
where $f$ and $P_l$ are as in the hypothesis of the theorem. On a certain family $V_j$ of sectors covering a disk punctured at $0\in\C$, we construct functions
$F_j:V_j\rightarrow\C$ having $\widehat{F}$ as their common asymptotic expansion.

We show two results. On one hand, if
we suppose that $\widehat{F}(x)$ is identically zero, then, by
Ramis-Sibuya's Lemma and the Relative Watson's Lemma,
any $F_j$ is exponentially flat of order strictly
greater than $p\nu$ where $p$ is the Poincar\'{e} rank of the
system and $\nu$ is the minimum of the valuations of the polynomials
$P_l$. On the other hand, assuming the condition (\ref{stokes})
(and a preliminary normalization of the function $f$), we show that
if $f$ is not identically zero,
then at least one of the differences $F_{j+1}-F_j$ is
exponentially flat of order exactly $p\nu$. Thus, by the first result, $f$ must vanish identically if $\widehat{F}(x)\equiv 0$.

\subsection{Examples and applications}
In the following examples and applications we discuss the
conditions (\ref{eq:argumento-autovalores}) and (\ref{stokes}) and
their relations with the properties (SQA) and (SAT).

\begin{example} {\em It is not sufficient in Theorem \ref{main3}
(and Theorem \ref{main2} as well) to assume only that $H(x)\rightarrow0$
as $x\rightarrow0$. In the example
$$x^2y'=\left[\matr{rr}{0&1\\-1&0}+\mu\, x\,id\right]y+\matr{r}{x\\x},\ \ \mu>0,$$
the eigenvalues of the leading matrix are $\pm i$ and hence
the unique formal solution is 1-summable in the positive real direction;
conditions (\ref{eq:argumento-autovalores}) and (\ref{stokes}) are also satisfied
(compare examples \ref{ex:autovalores-complejos}, \ref{ex:enlazamiento}).
Moreover, the system has a homogeneous solution of the form
$x^\mu(\sin(1/x),\cos(1/x))^T$ tending to 0 as $x\rightarrow 0,$
whose graph obviously cannot belong to any o-minimal structure
and which {\em does not have an asymptotic expansion} in a power series
as $x\rightarrow0.$ By superposition, the same is true for the above
non-homogeneous system.}
\end{example}\medskip

\begin{example}[Plane Pfaffian curves]\labe{ex:curvas-pfaffianas}
{\em A first consequence of our results is the model completeness
of the structures generated by certain pfaffian curves.
Consider a 1-form $\omega=a(x,y)dx+b(x,y)dy$ real analytic in some
neighborhood of $0\in\R^2$ with $a(0,0)=b(0,0)=0$.
A {\em pfaffian curve} defined by $\omega$ is the graph $\cal C$
of a non-oscillating solution $H:(0,\eps]\rightarrow\R$ tending
to 0 as $x\rightarrow 0$ of the differential equation
$b(x,y)\frac{dy}{dx}=-a(x,y)$\ .
The o-minimality of the structure $\R_{an,H}$ is a particular case of
the results proved in \cite{Wil,Lio-Rol}.

Analytic changes of coordinates or blowups are obviously
inessential for the question of model completeness.
Therefore, by a classical theorem on the reduction of singularities
\cite{Sei}, we can suppose that the origin is a simple singularity of
$\omega$ and hence that $H$ is of one of the following types:
\begin{enumerate}\item $H$ is analytic also at $x=0$. In this case, the model completeness of $\R_{an,H}$ is proved in \cite{Gab2}.
\item $H(x)=x^\lambda$ for some irrational $\lambda>0$. In this case,
the model completeness is proved in \cite{Mil}.
\item $H$ is a solution of some {\em saddle-node equation}
\begin{equation}\labe{eq1}
x^{p+1}\frac{dy}{dx}=y+A_1(x,y),\;\;p\geq 1,
\end{equation}
where $A_1$ is analytic at the origin and $A_1(0)=\frac{\partial A_1}{\partial y}(0)=0$. In this case, if the Taylor series $\widehat H(x)$ is divergent,
the conditions of Theorem \ref{main3} are satisfied, because a formal solution
$\widehat H(x)$ without a singular direction having non-trivial Stokes
phenomenon would be convergent (see also \cite{Mar-Ram}).
Thus the structure $\R_{an,H}$ is model complete.
If $\widehat H(x)$ converges, however, the model completeness is still an open question.
\end{enumerate}}
\end{example}

\begin{example} \labe{ex:euler2}

{\em Assume that (\ref{sistema}) has a formal power series solution $\widehat{H}(x)\in\R[[x]]^r$. Contrary to the scalar case, if $r>1$ then
divergence of $\widehat{H}(x)$ does not suffice for strong quasi-analyticity of solutions with asymptotic expansion $\widehat{H}(x)$. Even the
condition that $\widehat{H}(x)$ is {\em analytically transcendental}, i.e.\ $f(x,\widehat{H}(x))$ is not identically zero for all non zero analytic
function $f(x,\yy)$, does not imply strong quasi-analyticity (compare with (SAT)).

For example, let $E:x\mapsto E(x)$, $x>0$, be a solution of the
Euler equation $x^2\frac{dy}{dx}=y-x$ and define the function $H:(0,\infty)\rightarrow\R^2$ by
$H(x)=(E(x)+{\rm \exp}(-1/x),E(2x))$. It satisfies the system of
differential equations
\begin{equation}\labe{eq:euler2}
x^2 \frac{dy_1}{dx}=   y_1-x,\;\; \;\;
x^2 \frac{dy_2}{dx}= \frac{y_2}{2}-x.
\end{equation}
Clearly, $H$ does not satisfy condition (SQA). However, its asymptotic expansion $\widehat{H}(x)=(\widehat{E}(x),\widehat{E}(2x))$ diverges,
$\widehat{E}(x)=\sum_{n=1}^{\infty}(n-1)!x^{n}$ being the {\em Euler series}. Moreover, since $\widehat{E}(x)$ satisfies the property (SAT) by
Theorem \deuxbis{} applied to the Euler equation, $\widehat{H}(x)$ is analytically transcendental.

Notice that the linear part of (\ref{eq:euler2}) has eigenvalues
$1,1/2$, so that condition (\ref{eq:argumento-autovalores}) is not
satisfied.
}
\end{example}

\begin{example}[Two dimensional systems with non-real eigenvalues]\labe{ex:autovalores-complejos}
\mbox{\ }\\ {\em Consider a system (\ref{sistema}) with $r=2$ such that the linear
part $A_0=\partial A/\partial\yy(0)$ has two non real conjugate
eigenvalues $\lambda_1,\lambda_2=\overline{\lambda_1}$ (in
particular, it satisfies condition
(\ref{eq:argumento-autovalores})). Let
$\widehat{H}(x)\in\R[[x]]^2$ be its formal power solution. The
singular directions $d_{\theta_{l,j}}$ of $\widehat{H}(x)$ satisfy
$\theta_{l,1}=-\theta_{k,2}$ if $l\equiv -k$ ${\rm mod}\,p$ and,
since $\widehat{H}(x)$ has real coefficients, the corresponding
Stokes coefficients $c_{l,1}$, $c_{k,2}$ are complex conjugate.
Thus, condition (\ref{stokes}) is equivalent to the existence of a
non zero Stokes coefficient for some singular direction and hence
to the divergence of the series $\widehat{H}(x)$. We can summarize
the preceding discussion by stating the following
\begin{corollary}\labe{cor:autovalores-complejos}
Consider a two dimensional system (\ref{sistema}) such that the linear part has no real eigenvalues and suppose that its formal power series solution
$\widehat{H}(x)\in\R[[x]]^2$ does not converge. Then $\widehat{H}(x)$ is strongly analytically transcendental (SAT). As a consequence, any solution
$H:(0,\varepsilon]\rightarrow\R^2$ with $H(x)\sim\widehat{H}(x)$ as $\R^+\ni x\ra0$ is strongly quasi-analytic and the structure $\R_{an,H}$ is
o-minimal and model complete.
\end{corollary}
}
 \end{example}

\begin{example} {\bf (Non-oscillating solutions and o-minimal structures)}
\labe{ex:enlazamiento}\linebreak {\em Two dimensional systems (\ref{sistema}) with non real eigenvalues as in Corollary
\ref{cor:autovalores-complejos} were already treated in \cite{Can-Mou-San} in the context of analytic three dimensional vector fields. There, the
authors study the whole family $IP$ (``integral pencil") of solutions $H:(0,\varepsilon]\rightarrow\R^2$ of (\ref{sistema}) having $\widehat{H}(x)$
as the asymptotic expansion as $\R^+\ni x\ra0$. They show that any element in $IP$ is non-oscillating if and only if $\widehat{H}(x)$ is divergent.
Moreover, in that case, $\widehat{H}(x)$ is analytically transcendental and any pair of solutions in $IP$ are {\em asymptotically linked} in the
sense that each one ``turns around" the other.

Let us apply our results to an example
considered in \cite{Can-Mou-San}, generalizing
the Euler equation to dimension 2:
\begin{equation}\labe{eq2}
  \left\{
  \begin{array}{l}
   x^{2}\frac{dy_1}{dx}=y_1+y_2-x, \\
   x^{2}\frac{dy_2}{dx}=y_2-y_1.
  \end{array}
  \right.
 \end{equation}
 The eigenvalues of its linear part are $1\pm i$. It is easily seen
 that any solution
 $H(x)=(H_1(x),H_2(x))$ of (\ref{eq2}) defined for $x>0$ verifies
 $\lim_{x\rightarrow 0}H(x)=0$ and that all such solutions have a
 common asymptotic expansion
 $\widehat{H}(x)\in\R[[x]]^2$ as $\R^+\ni x\ra0$. Moreover, $\widehat{H}(x)$
 is divergent; this follows from the fact that $\hat y(x)=\widehat H_1(x)+i\,
 \widehat H_2(x)$ is the formal solution of a modification
 $x^2y'=(1-i)y-x$ of Euler's equation and hence
 $\hat y(x)= \sum_{n=1}^\infty (n-1)!\,(1-i)^{-n}\, x^n $ diverges.
By Corollary \ref{cor:autovalores-complejos}, the
structure $\R_{an,H}$ is o-minimal and model complete for any such
solution.

 The most important feature
 of this example is the following:
 given two {\em different} solutions $H,G:(0,\varepsilon)\rightarrow\R^2$ of (\ref{eq2}),
 the argument of the vector
 $x\mapsto H(x)-G(x)$ tends to infinity as $x\ra 0$
(``asymptotic linking"). As a consequence, the two structures $\R_{an,H}$,
 $\R_{an,G}$ do not admit any common o-minimal extension.
 Thus, we obtain an explicit
 family
 of ``mutually incompatible" o-minimal
 structures. Examples of incompatible pairs of o-minimal structures
were already exhibited in
 \cite{Rol-Spe-Wil}.

 It is natural to wonder if the o-minimality of these structures could arise
from well-known results of pfaffian geometry. It is proved in \cite{speiss:clos} that the so-called {\em pfaffian closure}
$\mathcal{S}_{\text{Pfaff}}$ of any o-minimal structure $\mathcal{S}$ is o-minimal as well. Because of the previous remark, it is clear that at most
one of the structures $\mathbb{R}_{an,H}$ could be a reduct of the pfaffian closure $\mathbb{R}_{an,\text{Pfaff}}$ of $\mathbb{R}_{an}$. We claim
that it is actually the case for none of these structures.

Let us suppose indeed that for some solution $H$ of example $( 2.5 )$ the structure $\mathbb{R}_{an,H}$ is a reduct of
$\mathbb{R}_{an,\text{Pfaff}}$. Then, for any other solution $G$, the pfaffian closure $\mathbb{R}_{an,G,\text{Pfaff}}$ (which is of course an
extension of $\mathbb{R}_{an,\text{Pfaff}}$) would be a common o-minimal extension of $\mathbb{R}_{an,H}$ and $\mathbb{R}_{an,G}$, which is
impossible.

 Example (\ref{eq2}) also provides an example of a non-oscillating solution of
 systems of analytic ordinary differential equations which is not definable in
any o-minimal extension of $\R$:

  \begin{proposition}\labe{contraejemplo}
 Let $H, G:(0,\varepsilon)\rightarrow\R^2$ be two different solutions of (\ref{eq2}).
 Then, the function $H^*(x)=(H(x),G(2x)):(0,\varepsilon)\rightarrow\R^4$ is a non-oscillating
 solution of
 a system of analytic ordinary differential equations of the form (\ref{sistema}) such that
 the structure $\R_{H^*}$ is not
 o-minimal.
  \end{proposition}
 \begin{proof} The asymptotic expansion
 $\widehat{H}^*(x)=(\widehat{H}(x),\widehat{H}(2x))$ of $H^*(x)$ as $\R^+\ni x\ra0$
 is analytically transcendental by Theorem \deuxbis{}. This implies that $H^*(x)$
 is non-oscillating. On the other hand, the image of the curve $x\mapsto
 (H(x),G(x))$ in $\R^4$ is definable in the structure $\R_{H^*}$ but it cuts
 infinitely many times the semialgebraic set $\{(v_1,v_2)\in\R^4/v_1-v_2=0\}$.
\end{proof}

 We finish with some application of Proposition \ref{contraejemplo}
 to {\em Hardy fields}. The latter are, by definition,
 subfields of the ring of germs of differentiable real functions on semi-intervals
 $(0,\varepsilon)$
 which are closed under the usual differentiation (see \cite{Ros}).
 A non-oscillating solution $H(x)$ of a system (\ref{sistema}) determines a
 Hardy field $K_H$ consisting of the germs of quotients
 $f(x,H(x))/g(x,H(x))$ where $f,g$ are analytic functions and $g(x,H(x))\not\equiv
 0$. Another example of a  Hardy field is
 the field $K_\RR$ of germs of functions of one variable
definable in some o-minimal structure
$\RR$ \cite{Lio-Mil-Spe}.
Proposition~\ref{contraejemplo} gives an example
of a {\em polynomially bounded} Hardy field which
is not contained in the field $K_\RR$ for any o-minimal structure $\RR$. An example
of this kind, but less explicit, is already exhibited in
\cite{Bos}.
}\end{example}

\subsection{Open questions}

We have three main questions:

1) Given a solution $H:(0,\varepsilon)\rightarrow\R$ of a scalar equation as (\ref{eq1}) with a convergent asymptotic expansion at $x=0$, is the
structure $\R_{an,H}$ model complete? We conjecture that this is not always the case.

2) Does any non-oscillating solution
$H:(0,\varepsilon)\rightarrow\R^r$ of a system of analytic ordinary differential equations
(\ref{sistema}) with $r=2$ or $r=3$ generate an o-minimal
structure?

3) Which are the relations between analytic transcendence, quasianalyticity, the hypotheses (\ref{eq:argumento-autovalores}) and (\ref{stokes}) and
(non linear) differential Galois theory?

\section{Quasi-analyticity and model-completeness}
\setcounter{theorem}{0}\setcounter{equation}{0}

This paragraph is devoted to the proof of Theorem \ref{main1}.

\subsection{Some generalities}

Let $H=(H_1,\ldots,H_r):(0,\varepsilon]\rightarrow\R^r$ be a solution of a
system of (real) analytic differential equations of the form (\ref{sistema}) such that
$\lim_{x\rightarrow 0}H(x)=0$. The following general result, which will be
essential in the sequel, might already be known. However, we have not
found any reference.
\begin{lemma}\labe{morfismo}
Suppose that $A(0,0)=0$. Then, for any $L>0$ there exists a neighborhood $V$ of
$0\in\R^{r+1}$, $\delta_L>0$ and an analytic function $B:[-L,L]\times
V\rightarrow\R^{r}$ such that
$$
H(x+x^{p+1}z)=B(z,x,H(x)),\;\;\; \mbox{for }| z|\leq L\mbox{ and
}0<x\leq\delta_L.
$$
\end{lemma}
Observe that, conversely, any function $H$ satisfying the above functional
equation also solves a certain system (\ref{sistema}) of differential equations,
namely $x^{p+1} H'= \frac{\partial B}{\partial z}(0,x,H)$.

\begin{proof}
Let $V_1$ be a neighborhood of the origin in $\R^{r+1}$ such that the function
$A(x+x^{p+1}z,{\bf w})$ is defined for $(z,x,{\bf w})\in(-L-1,L+1)\times V_1$ and $1+x^pz$ does not vanish there.
Consider the regular analytic differential equation with analytic parameter
$x$:
$$
(E_x)\;\;\;\;\frac{d{\bf w}}{dz}=(1+x^pz)^{-p-1}A(x+x^{p+1}z,{\bf w}).
$$
Since $A(0,0)=0$, the solution of the initial value problem ($E_0$), $w(0)=0$
is the zero function $w:(-L-1,L+1)\rightarrow\R$, $z\mapsto 0$. By the
theorem of analytic dependence upon parameters and initial values, there exists
a neighborhood $0\in V\subset V_1$ such that for $(x,w_0)\in V$ the solution of
($E_x$), $w(0)=w_0$ exists on the interval $[-L,L]$. Moreover, the mapping
$B:[-L,L]\times V\rightarrow\R^{r}$ associating $(z,x,w_0)$ to $w(z)$, the
value of this solution at some $z\in[-L,L]$, is analytic. Choose finally
$\delta_L>0$ such that $(x,H(x))\in V$ for $0<x\leq\delta_L$. Then $G_x(z)=B(z,x,H(x))$
is defined for $z\in[-L,L],\,0<x<\delta_L$ and denotes the solution of $(E_x)$,
$G_x(0)=H(x)$ which is thus defined on the interval $[-L,L]$.
For sufficiently small $z$,
this solution can also be expressed by $G_x(z)=H(x+x^{p+1}z)$. The identity theorem
for (real) analytic functions yields the theorem.
\end{proof}

Suppose now that $H$ has an asymptotic expansion
$$H(x)\sim\widehat{H}(x)=(\widehat{H}_1(x),\ldots,\widehat{H}_r(x))
        \in\R[[x]]^r$$
as $\R^+\ni x\ra0$.
\begin{lemma}\labe{lema:h-cinfinito}
The  function $H$ can be continued to a $C^\infty$ function on $[0,\varepsilon]$
and its Taylor series at $x=0$ is equal to
$\widehat{H}(x)$.
\end{lemma}
\begin{proof}
The function $x\mapsto A(x,H(x))$, $x>0$ has  $A(x,\widehat{H}(x))$ as its
asymptotic expansion as $0<x\ra0$. Hence, $x\mapsto
H'(x)=A(x,H(x))/x^{p+1}$ admits the asymptotic expansion
$A(x,\widehat{H}(x))/x^{p+1}$, which might contain
negative powers of $x$. Since asymptotic expansions are compatible
with integration, this cannot be the case and
$A(x,\widehat{H}(x))$ is divisible by $x^{p+1}$.
Moreover, we have the following equation for formal series:
$\widehat{H}(x)=\int_{0}^{x}A(t,\widehat{H}(t))t^{-p-1}dt$ and hence
$\widehat{H}(x)$ (formally)  satisfies the system (\ref{sistema}).
This implies that $H'(x)\sim\widehat{H}'(x)$ as $\R^+\ni x\ra0$. We conclude
that $H$ is of class $\CC^1$ at $x=0$ and that $\lim_{x\rightarrow
0^+}H'(x)=H'(0)=\widehat{H}'(0)$. The same argument for
$H'(x)$, which is a solution of $x^{p+1}z'=\frac{\partial A}{\partial x}(x,H(x))+
\frac{\partial A}{\partial y}(x,H(x))z-(p+1)x^pz$, shows that $H$ is of class $\CC^2$ at $x=0$. The
proof can be completed by induction in this way.
\end{proof}
\subsection{Simple functions and the class $\AAA$}
Lemma \ref{lema:h-cinfinito} allows to reduce the statement of Theorem
\ref{main1} to $\CC^\infty$ solutions defined on intervals symmetric to 0.
Indeed, consider
$G:[-\varepsilon,\varepsilon]\rightarrow\R^r$ defined by
$G(x)=H(x^2)$ if $x\neq 0$ and $G(0)=0$. This is a $\CC^\infty$ function
by Lemma \ref{lema:h-cinfinito} and it satisfies a system of ordinary differential equations
of the form (\ref{sistema}). Denote by
$\widetilde{H}$, $\widetilde{G}:\R\mapsto\R^r$, respectively, the extensions of $H$
resp.\ $G$ that vanish outside their corresponding domains of definition.
In a straightforward manner, it can be shown
that $\R_{\widetilde{F}_{an}\cup\{\widetilde{G}\}}$ is o-minimal
(respectively model complete) if and only if
$\R_{\widetilde{F}_{an}\cup\{\widetilde{H}\}}$ is o-minimal
(respectively model complete).

Moreover, if $H$ satisfies the strong quasi-analyticity property (SQA) then
$G$ has the analogous property without the condition that
the first non-zero coefficients $P_l^{({\rm val} P_l)}(0)$ of the polynomials $P_l$ have to be positive.

We introduce some notation useful for the proofs (recall that
$T_k\phi\,(x)=(\phi(x)-J_k\phi\,(x))/x^{k}$ for a $\CC^\infty$ function
or a formal series $\phi$, where $J_k\phi\,(x)\in\R[x]$ denotes its $k$-jet at 0).
\begin{definition}
Let $H=(H_1,\ldots,H_r):[-\varepsilon,\varepsilon]\rightarrow\R^r$ be a
$\CC^\infty$ function. A germ $\varphi(x)$ of a $\CC^\infty$ function at $0\in\R$
will be called a {\em simple function} (relatively to $H$) if there exists
$n\geq 0$, an analytic function $f\in\R\{x,z_{11},\ldots,z_{rn}\}$, polynomials
$P_j(x)$ with $P_j(0)=0$ for $j=1,\ldots,n$ and an integer $k\geq 0$ such that
\begin{equation}\labe{eq-simples}
\varphi(x)=f(x,T_k H(P_1(x)),\ldots,T_k H(P_n(x))).
\end{equation}
The family of simple functions is an algebra denoted by $\SC_H$ (or simply by
$\SC$ if $H$ is obvious from the context).
\end{definition}

Using the above considerations, Theorem \ref{main1} follows immediately
from the following seemingly weaker statement:
\begin{theorem}\labe{main11}
Let $H=(H_1,\ldots,H_r):[-\varepsilon,\varepsilon]\rightarrow\R^r$
be a $\CC^\infty$ solution of a system  (\ref{sistema}) of ordinary differential equations.
Suppose that the algebra $\SC$ of simple functions relatively to $H$ is
quasi-analytic; i.e.\ $\varphi=0$ whenever the Taylor series
of some $\varphi\in\SC_H$ at $x=0$ vanishes.
Then the extension $\R_{an,H}$ is o-minimal and model complete.
\end{theorem}

Our proof of Theorem  \ref{main11}, which will be given in subsection 3.4,
relies on a fundamental result of \cite{Rol-Spe-Wil}
about the model completeness of structures generated by certain
quasi-analytic classes. Therefore we first study the quasi-analyticity of
certain classes of functions defined using $H$.

Given a $\CC^\infty$ function
$H=(H_1,\ldots,H_r):[-\varepsilon,\varepsilon]\rightarrow\R^r$,
let $\AAA_H$ be the smallest collection  $\AAA_{H}^{m}$, $m\in\N$, of
subalgebras of germs of $\CC^\infty$ functions at
$0\in\R^m$ satisfying the following conditions
(we drop the subscript ``$H$" for the sake of clarity):
\begin{itemize}
    \item[A1)] The germs of analytic functions 
    of $m$ variables are elements of $\AAA^{m}$,
    for all $m$, and the germ of $H_l$ at $0$ is in $\AAA^{1}$
    for
    $l=1,\ldots,r$.
    \item[A2)] (Stability by composition). If $\phi\in\AAA^{m}$
    and $\phi_1,\ldots,\phi_m\in\AAA^{n}$ with
    $\phi_1(0)=\cdots=\phi_m(0)=0$ then
    $\phi(\phi_1,\ldots,\phi_m)\in\AAA^{n}$.
    \item[A3)] (Stability by implicit equations). If $\phi\in\AAA^{m+1}$
    satisfies $\phi(0)=0$ and $\partial\phi/\partial x_{m+1}\,(0)\neq 0$
    then there exists $\bar{\phi}\in\AAA^{m}$ such that
    $\bar{\phi}(0)=0$ and the germ $\phi(\xx,\bar{\phi}(\xx))\in\AAA^{m}$,
    with $\xx=(x_1,\ldots,x_m)$, is
    equal to zero.
    \item[A4)] (Stability by monomial division). If $\phi\in\AAA^{m}$
    satisfies $\phi(0,\bar{\xx})\equiv 0$, where $\bar{\xx}=(x_2,\ldots,x_m)\in\R^{m-1}$
    then there exists $\bar{\phi}\in\AAA^{m}$ such that $\phi(\xx)=x_1
    \bar{\phi}(\xx)$,
    where $\xx=(x_1,\bar{\xx})$.
\end{itemize}
\begin{remark}{\em Properties A2-A4 imply that the algebras
$\AAA^{m}$ are stable by partial differentiation. In fact, if $\phi\in\AAA^{m}$, the
germ $\phi_1(u,\xx)=\phi(x_1+u,\bar\xx)-\phi(\xx)$ is in $\AAA^{m+1}$ (here again
$\xx=(x_1,\bar{\xx})$), and thus also
$\bar{\phi}_1(u,\xx)=\phi_1(u,\xx)/u$ defines an element of $\AAA^{m+1}$ by A4. Thus,
$\partial\phi/\partial x_1(\xx)=\bar{\phi}_1(0,\xx)$ is in $\AAA^{m}$.}
\end{remark}

\begin{remark}{\em Property A3 of the class $\AAA$ can be generalized to systems
of implicit equations:

{\em
Let $\Phi=(\Phi_1,\ldots,\Phi_n)\in(\AAA^{m+n})^n$ such that
$\Phi(0)=0$ and $\partial\Phi/\partial\yy(0)$ is an invertible
matrix (with the notation $(\xx,\yy)$ of the coordinates of
$\R^{m+n}$). Then there exists
 $\phi=(\phi_1,\ldots,\phi_n)\in(\AAA^m)^n$
such that $\Phi(\xx,\phi(\xx))\equiv 0$.
}

\begin{proof}
 By induction over $n$. The case $n=1$ is property A3. So suppose that $n>1$.
  Denote $\partial \Phi/\partial\yy(0)= (a_{ij})_{1\leq i,j\leq n}$.
 Up to some permutation of the rows, we can suppose that its principal minors
 ${\rm det}\,((a_{ij})_{1\leq i,j\leq k})$, $k=1,\ldots,n$, are all
 different from zero. By property A3, the solution $\psi_1$ of
 $$
 \Phi_1(\xx,\psi_1(\xx,y_2,\ldots,y_n),y_2,\ldots,y_n)=0
 $$
 is in $\AAA^{m+n-1}$. For $j=2,\ldots,n$ consider
 $\bar{\Phi}_j(\xx,\bar{\yy})=\Phi_j(\xx,\psi_1(\xx,\bar{\yy}),\bar{\yy})$ and
 $\bar{\Phi}=(\bar{\Phi}_2,\ldots,\bar{\Phi}_n)$, where
 $\bar{\yy}=(y_2,\ldots,y_n)$. We have that $\bar{\Phi}_j\in\AAA^{m+n-1}$ and
 $$
\frac{\partial \bar{\Phi}_j}{\partial y_k}(0)=\frac{1}{a_{11}}{\rm det}\,
\left(
\begin{array}{cc}
  a_{11} & a_{1k} \\
  a_{j1} & a_{jk} \\
\end{array}\right),
 $$
 so that $\bar{J}=\partial \bar{\Phi}/\partial\bar{\yy} (0)$ is invertible. By the hypothesis of
 induction, the system $\bar{\Phi}(\xx,\bar{\phi}(\xx))=0$ is solved by
some $\bar{\phi}=(\phi_2,\ldots,\phi_n)\in(\AAA^{m})^{n-1}$ and completing the vector
with $\phi_1(\xx)=\psi_1(\xx,\bar{\phi}(\xx))\in\AAA^{m}$ we have solved
$\Phi(\xx,\phi(\xx))\equiv 0$.
\end{proof}
}
\end{remark}
We pass from germs to actual functions by means of the following definition. A
function $F:U\rightarrow\R$ on an open subset $U\subset\R^m$ is called an {\em
$\AAA$-analytic function} if for any $a\in U$ there exists $\phi_a\in\AAA^{m}$
such that the germ of $F$ at $a$ is equal to the germ of
$\xx\mapsto\phi_a(\xx-a)$. This definition is justified by the following:
\begin{lemma}\labe{lema-germenes}
If $\phi\in\AAA^{m}$, there exists a representative of $\phi$ defined in some
neighborhood of $0\in\R^m$ that is an $\AAA$-analytic function.
\end{lemma}
\begin{proof}
For any $m$, let $\widetilde{\AAA}^{m}$ be the subset of $\AAA^{m}$ consisting
of all germs having the wanted property. One can verify easily that the
collection $\widetilde{\AAA}=(\widetilde{\AAA}^{m})_m$ verifies properties
A1--A4. Thus, $\widetilde{\AAA}^{m}=\AAA^{m}$ for any $m$.
\end{proof}

\subsection{Quasi-analyticity of the class $\AAA$}
Clearly, if $\AAA=(\AAA^m)$ is a family of algebras of germs of $\CC^\infty$
functions that satisfies properties A1--A4 then the collection of simple
functions with respect to $H$ is a subalgebra of $\AAA^1$. Our goal in this
paragraph is to prove that if $H$ is a solution of a system of analytic
ordinary differential equations,
then the quasi-analyticity property of $\SC_H$ is inherited by the entire class
$\AAA_H$.
\begin{theorem}[Quasi-analyticity of $\AAA$] \labe{QA}
Under the hypothesis of Theorem \ref{main11}, let $\AAA=\AAA_H$ be
the smallest class of germs satisfying A1--A4. Then $\AAA$ is a
quasi-analytic class: for any $m\geq 1$, if the Taylor series
of $\phi\in\AAA^m$ at $0$
vanishes then $\phi=0$.
\end{theorem}

The proof is given in the remainder of this subsection and it is
divided in two parts. In the first one, we deduce
the general statement from the one for
germs of a single variable. In the
second part, we show the surprising fact that any element of $\AAA^1$ is a simple
function.
\begin{proposition}\labe{prop:1-variable}
If the algebra $\AAA^1$ of elements of $\AAA$ of one variable is
quasi-analytic then Theorem \ref{QA} follows.
\end{proposition}
\begin{proof}
Let $\phi\in\AAA^m$ such that its Taylor series
$\hat{\phi}(\xx)\in\R[[\xx]]$ vanishes. By Lemma
\ref{lema-germenes}, we can consider a representative $F$ of
$\phi$ which is an $\AAA$-analytic function in the open ball
$B(0,\epsilon)$, $\epsilon>0$. Given  a unitary vector $a\in\R^m$,
let $F_a:(-\epsilon,\epsilon)\rightarrow\R$ be defined by
$F_a(t)=F(ta)$. Then the germ of $F_a$ at $t=0$ is an element of
$\AAA^1$ whose Taylor series vanishes. Our assumption implies
that the function $F_a$ vanishes identically in some interval
centered at $t=0$. Let $\delta_a=\sup\{\delta\in
(0,\epsilon)/F_a(t)=0\mbox{ for }0\leq t\leq\delta\}>0$. The proof
is complete if we show that $\delta_a=\epsilon$ for any $a$.
Suppose, on the contrary, that $\delta_a<\epsilon$ for a certain
$a$. Then the function $G_a$ defined by $G_a(t)=F_a(t+\delta_a)$,
defined on $]-\delta_a,\epsilon-\delta_a[$, vanishes for $t<0$. Thus, its
Taylor series at $t=0$ must vanish. Since the germ of $G_a$ at
$t=0$ belongs to $\AAA^1$ by choice of $F$, the function $G_a$ vanishes in a
neighborhood of $0$. This contradicts the definition of
$\delta_a$.
\end{proof}

\begin{proposition}\labe{prop:a1=simples}
Under the hypothesis of Theorem \ref{main11}, we have $\AAA^1=\SC_H$.
\end{proposition}
As $\SC_H$ is assumed to be quasi-analytic in the hypothesis of Theorem
\ref{main11} (which is also that of Theorem \ref{QA}), propositions
\ref{prop:1-variable} and \ref{prop:a1=simples} prove Theorem  \ref{QA}.
It remains to prove the above proposition.

This result is by no means trivial and actually quite surprising: the
definition of the class $\AAA$ permits solutions of implicit
equations and compositions, so that we could expect that $\AAA^1$ contains
much more complicated functions than the simple ones. The proof of
Proposition \ref{prop:a1=simples} needs several preparatory lemmas
and will be given at the end of this subsection. The hypothesis of
Theorem \ref{main11} (i.e.\ the quasi-analyticity of the class
$\SC_H$ of simple functions)
will be assumed in the remainder of this subsection
without explicitly mentioning it.
\begin{lemma}\labe{lema:simples}
The class of simple functions $\SC_H$ is closed by monomial
division and by composition.
\end{lemma}
\begin{proof}
Let $\varphi\in\SC_H$ be written as
$\varphi(x)=f(x,T_kH\,(P_1(x)),\ldots,T_kH\,(P_n(x)))$ with
$f$ analytic and polynomials $P_j(x)$ such that $P_j(0)=0$.
Suppose that $\varphi(0)=0$. Using the formula $T_kH\,(x)=x{\bf
a}_{k+1}+xT_{k+1}H\,(x)$, where ${\bf a}_{k+1}\in\R^r$ is the
$k+1$-th coefficient of the Taylor series $\widehat{H}(x)$ of
$H$ at $0$, we can write the quotient $\varphi(x)/x$ in the form
$$
\tfrac{1}{x}f(x,\ \{P_l(x){\bf a}_{k+1}+P_l(x)T_{k+1}H\,(P_l(x))\}_{l=1,...n}).
$$
This last function is of the form
$\tilde{f}(x,T_{k+1}H\,(P_1(x)),\ldots,T_{k+1}H\,(P_n(x)))$ with
some analytic $\tilde{f}$. This proves the closure of the simple functions by
monomial division.

\noindent Consider now $\varphi_1,\varphi_2\in\SC_H$ such that
$\varphi_2(0)=0$ and let us prove that the composition
$\varphi=\varphi_1\circ\varphi_2$ is also a simple function. It
suffices to consider the case that $\varphi_1$ is a component of $H$, say $H_1$.
Assume that $\varphi_2\neq 0$ (otherwise there is nothing to
prove). Since $\SC_H$ is a quasi-analytic class, there exists
$d>0$ and $c_d\neq 0$ such that we can write
$$
\varphi_2(x)=c_dx^d+\cdots+x^{(p+1)d}\,T_{(p+1)d}\,\varphi_2\,(x)=
P(x)+P(x)^{p+1}\widetilde{\varphi}_2(x),
$$
where
$P(x)=c_dx^d+\cdots+c_{pd+d-1}x^{pd+d-1}=J_{pd+d-1}\,\varphi_2$
and $\widetilde{\varphi}_2$ is still a simple function as we already
proved that $\SC_H$ is closed under monomial division. By Lemma
\ref{morfismo}, we have
$\varphi(x)=H_1(\varphi_2(x))=B_1(x,H_1(P(x)),\widetilde{\varphi}_2(x))$
with some analytic function $B_1$. We conclude that
$\varphi\in\SC_H$.
\end{proof}

It is often convenient to view
a function $f$ obtained by the composition of several functions
$f_i$ in the following way: consider the graph of $f$ as the
projection of a set essentially constructed
using the graphs of the $f_i$ and cartesian products
(see for example \cite{Mou-Roc}).

The elements of $\AAA$ are obtained from $H$ and analytic functions
by repeated compositions, monomial divisions and solutions of implicit
equations. In the same spirit, in order to replace the recursive definition
of elements of $\AAA$, we describe them as projections of
solutions of (big) systems of (simple) implicit equations in
several variables. The individual equations can be analytic or in a new class of functions
slightly extending $H$, the {\em divided differences}.

We define a divided difference
(generated by $H$) as an element of the smallest subcollection
$\DD=(\DD^m)$ of $\AAA$ containing the germ of each component
$H_l$ of $H$ and such that, if $\Delta$ is in $\DD^m$,
then the germ $\bar{\Delta}$ defined by one of the following
formulas also is in $\DD$:
\begin{itemize}
   \item[D1)] $\bar{\Delta}(\xx)=\Delta(0,x_2,\ldots,x_m)$.
   \item[D2)] $\bar{\Delta}(\xx)=\Delta(\xx)-\Delta(0)$ or $\bar{\Delta}(\xx)=\Delta(\xx)-\Delta(0,x_2,\ldots,x_m)$.
   \item[D3)] If $\Delta(0,x_2,\ldots,x_m)\equiv 0$ then
    $\bar{\Delta}(\xx)=\Delta(\xx)/x_1$.
   \item[D4)] $\bar{\Delta}(\xx)=\Delta(x_{\sigma(1)},\ldots,x_{\sigma(m)})$
    for some permutation $\sigma$ of $\{1,\ldots,m\}$.
   \item[D5)] $\bar{\Delta}(\xx,z)=\Delta(\xx)$.
  \item[D6)] $\bar{\Delta}(\xx,z)=\Delta(x_1+x_2(\alpha+z),x_2,\ldots,x_m)$
    for some $\alpha\in\R$.
\end{itemize}
\begin{lemma}\labe{lema2}
Let $\beta\in\AAA^m$, $m\geq 1$. Then there exists $n\in\N$ and a system of $n$
implicit equations
$\Phi(\xx,\yy)=(\Phi_1(\xx,\yy),\ldots,\Phi_n(\xx,\yy))^T=0$,
$(\xx,\yy)\in\R^{m+n}$ with $\Phi(0)=0$ and $det(\partial
\Phi/\partial\yy\,(0))\neq 0$ such that its solution
$\yy=\phi(\xx)=(\phi_1(\xx),\ldots,\phi_n(\xx))$ satisfies
$\phi_n=\beta-\beta(0)$ and such that each component $\Phi_i$ is either
analytic or of the form $\Phi_i(\xx,\yy)=y_{l_i}-\Delta_i(\xx,\yy)$ with some
$1\leq l_i\leq n$ and some divided difference  $\Delta_i$ generated by $H$.
\end{lemma}
\begin{proof}
Consider, for any $m$, the subset $\widetilde{\AAA}^m\subset\AAA^m$ of all germs
satisfying the statement of the lemma. Let us show that the collection
$\widetilde{\AAA}=(\widetilde{\AAA}^m)$ verifies properties A1--A4 of the
definition of the class $\AAA$.

Obviously, $\widetilde{\AAA}$ verifies A1.

Let $\beta_0\in\widetilde{\AAA}^m$ and
$\beta_1,\ldots,\beta_m\in\widetilde{\AAA}^n$ with $\beta_j(0)=0,\,j=0,...,m$.
We want to show that the composition
$\beta=\beta_0(\beta_1,\ldots,\beta_m)\in\widetilde{\AAA}^n$.

Consider the system
$\Phi^{(0)}(\zz,\yy^{(0)})=0$ of $n_0$ equations in $\R^{m+n_0}$ such that
its solution
$\yy^{(0)}=\phi^{(0)}(\zz)=(\phi^{(0)}_{1}(\zz),\ldots,\phi^{(0)}_{n_0}(\zz))$
defines $\beta_0=\phi^{(0)}_{n_0}$ and
consider, for $j=1,\ldots,m$, the system
$\Phi^{(j)}(\xx,\yy^{(j)})=0$  of $n_j$ equations in $\R^{n+n_j}$ such that its
solution
$\yy^{(j)}=\phi^{(j)}(\xx)=(\phi^{(j)}_{1}(\xx),\ldots,\phi^{(j)}_{n_j}(\xx))$
defines $\beta_j=\phi^{(j)}_{n_j}$. Put
$$
\Phi(\xx,\yy)=(\Phi^{(1)}(\xx,\yy^{(1)}),\ldots,\Phi^{(m)}(\xx,\yy^{(m)}),
\Phi^{(0)}((y^{(1)}_{n_1},\ldots,y^{(m)}_{n_m}),\yy^{(0)})),
$$
with
$\yy=(\yy^{(1)},\ldots,\yy^{(m)},\yy^{(0)})\in\R^{n_0+\cdots+n_m}$;
here the column vectors have to be stacked on top of each other.
We have $\Phi(0)=0$ and the matrix
$\partial \Phi/\partial\yy\,(0)$ is lower block diagonal with
diagonal blocks given by the matrices $\partial
\Phi^{(j)}/\partial\yy^{(j)}\,(0)$, thus invertible. If
$\yy=\phi(\xx)$ is the solution of the system $\Phi(\xx,\yy)=0$ then
$\phi=(\phi^{(1)},\ldots,\phi^{(m)},\psi)$ with $\psi$ satisfying
$\Phi^{(0)}((\phi^{(1)}_{n_1}(\xx),\ldots,\phi^{(m)}_{n_m}(\xx)),\psi(\xx))\equiv
0$. Hence, the last component of $\psi$ is
$\psi_{n_0}=\beta_0(\beta_1,\ldots,\beta_m)$. On the other hand,
it is obvious that each component $\Phi_i$ of $\Phi$ is either
analytic or of the form $y^{(j(i))}_{l_{i}}-\Delta_i(\xx,\yy)$ with
some divided difference $\Delta_i$. This proves A2 for
$\widetilde{\AAA}$.

Let $\beta(\xx,x_{m+1})\in\widetilde{\AAA}^{m+1}$ with
$\beta(0)=0$, $\partial \beta/\partial x_{m+1}\,(0)\neq 0$ and
consider the solution $\beta'(\xx)$ of
$\beta(\xx,\beta'(\xx))\equiv 0$. Let $\Phi((\xx,x_{m+1}),\yy)=0$
be the system of implicit equations corresponding to $\beta$ and
consider its solution $\yy=\phi=(\phi_1,\ldots,\phi_n)$ which satisfies
$\phi_n=\beta$. Put $\yy'=(\yy,y_{n+1})$ with an additional variable
$y_{n+1}$ and let
$$
\Phi'(\xx,\yy')=(\Phi((\xx,y_{n+1}),\yy)^T,y_{n})^T.
$$
We verify that $\Phi'(0)=0$ and, using $\partial \beta/\partial x_{m+1}\,(0)\neq 0$,
that the matrix
$$
\frac{\partial\Phi'}{\partial\yy'}(0)= \left(
\begin{array}{cc}
  \frac{\partial \Phi}{\partial\yy}(0) & {\bf b} \\
  {\bf e}_n & 0 \\
\end{array}\right),\;\;{\bf
b}=\partial\Phi/\partial x_{m+1}(0),{\bf
e}_n=(0,\ldots,0,1)\in\R^n
$$
is invertible. If $\yy'=\phi'=(\phi'_1,\ldots,\phi'_{n+1})$ is
the solution of the system $\Phi'(\xx,\yy')=0$ then we have
$\Phi((\xx,\phi'_{n+1}(\xx)),\,(\phi'_1(\xx),...\phi'_n(\xx))\equiv 0$ and
$\phi'_{n}(\xx)\equiv 0$.
The first equation implies that $\phi'_l(\xx)=\phi_l(\xx,\phi'_{n+1}(\xx))$
for $l=1,...,n$. The second equation then yields
$$\beta(\xx,\phi'_{n+1}(\xx))=\phi_n(\xx,\phi'_{n+1}(\xx))\equiv
0
$$
and thus, by definition, $\phi'_{n+1}=\beta'$. This proves A3 for
$\widetilde{\AAA}$.

Let $\beta\in\widetilde{\AAA}^m$ such that
$\beta(0,\bar{\xx})\equiv 0$ with $\bar{\xx}=(x_2,\ldots,x_{m})$
and let $\beta'\in\AAA^m$ be given by $\beta'(\xx)=\beta(\xx)/x_1\in\AAA^m$. Let
$\Phi(\xx,\yy)=0$ be a system of $n$ equations for $\beta$ given
by the lemma and $\yy=\phi=(\phi_1,\ldots,\phi_n)$ its solution.
Put $\psi=(\psi_1,\ldots,\psi_n)$ with
$\psi_j(\xx)=(\phi_j(\xx)-\phi_j(0,\bar{\xx}))/x_1\,\in\AAA^m$.
Consider the vector $ {\bf b}=\psi(0)=(\partial\phi/\partial
x_1\,(0))^T=-\frac{\partial \Phi}{\partial\yy}(0)\,^{-1}(\partial
\Phi/\partial x_1\,(0))$ and the system of $2n$ equations
$$
\Phi'(\xx,\yy,\yy')=(\Phi((0,\bar{\xx}),\yy),\tfrac{1}{x_1}[\Phi(\xx,\yy+x_1({\bf
b}+\yy'))-\Phi((0,\bar{\xx}),\yy)])^T.
$$
Denoting $\Psi(\xx,\yy,\yy')=\Phi(\xx,\yy+x_1({\bf b}+\yy'))$, we have
 $\Phi'(0)=(\Phi(0),\partial\Psi/\partial x_1\,(0))^T$ and this is the zero vector
 by the definition of ${\bf b}$. Furthermore, the matrix
 $\partial\Phi'/\partial(\yy,\yy')\,(0)$ is lower block triangular
 and its two diagonal blocks are equal to $\partial \Phi/\partial\yy\,(0)$,
 hence invertible. Moreover, each component of $\Phi$ is either analytic
 or a divided difference generated by $H$ in the variables
 $(\xx,\yy,\yy')$.
(Observe that this last part makes the introduction of the divided differences
unavoidable.)

We check that
 $\phi'(\xx)=(\phi(0,\bar{\xx}),\psi(\xx)-{\bf b})$ satisfies
 $\Phi'(\xx,\phi'(\xx))\equiv 0$. Its last component is
 $$
 \psi_n(\xx)-b_n=\tfrac{1}{x_1}(\phi_n(\xx)-\phi_n(0,\bar{\xx}))-b_n=
 \tfrac{1}{x_1}(\beta(\xx)-\beta(0,\bar{\xx}))-b_n=\beta'(\xx)-\beta'(0),
$$
and the proof is complete.
\end{proof}

In the following three lemmas, we establish some useful relations
between the divided differences and the simple functions.

\begin{lemma}\labe{lema3}
Let $\Delta(\xx)\in\DD^m$ be a divided difference in the variables $\xx=\xunoxm$.
Then there exists natural numbers $n,N\geq 0$, polynomials
$P_1(\xx),\ldots,$ $P_n(\xx),Q(\xx)$ with $Q(\xx)\not\equiv 0$ and a function
$B$ analytic in a neighborhood of $0\in\R^{m+nr(N+1)}$ such that
\begin{equation}\labe{eq-dd}
Q(\xx)\Delta(\xx)=B\left(\xx,\{H^{(l)}(P_j(\xx))\}_{0\leq l\leq N,1\leq j\leq
n}\right),
\end{equation}
where $H^{(l)}$ denotes the $l$-th derivative of $H$ with respect to $x$.
\end{lemma}
\begin{proof}
The components of $H$ trivially satisfy (\ref{eq-dd}). Moreover,
if $\Delta$ is a divided difference which satisfies (\ref{eq-dd}),
then the divided difference $\bar{\Delta}$ obtained from $\Delta$
by one of the operations D3-D6 also does. Also, if (\ref{eq-dd})
is true for
 $\Delta_1$ and $\Delta_2$ then it is true for the difference $\Delta_1-\Delta_2$.
  Let us now consider $\bar{\Delta}(\bar{\xx})=\Delta(0,\bar{\xx})$, where
$\bar{\xx}=(0,x_2,\ldots,x_m)$. Here the problem is that $Q(0,\bar{\xx})$ might
vanish identically. In any case, since $Q(\xx)\not\equiv 0$, we can consider the
smallest integer $s\geq 0$ such that $\partial^s Q/\partial
x_{1}^{s}\,(0,\bar{\xx})\not\equiv 0$. We differentiate the identity
(\ref{eq-dd}) $s$ times and obtain some
analytic function $\bar{B}$ in $m-1+nr(N+s+1)$
variables such that
$$
\frac{\partial^s Q}{\partial x_{1}^{s}}(0,\bar{\xx})\Delta(0,\bar{\xx})=
\bar{B}\left(\bar{\xx},\{H^{(l)}(P_j(0,\bar{\xx}))\}_{0\leq l\leq N+s,1\leq
j\leq n}\right).
$$
\end{proof}
\begin{lemma}\labe{lema4}
Let $\Delta\in\DD^m$ be a divided difference and
$\varphi_1,\ldots,\varphi_m\in\SC_H$ be simple functions. Then the composition
$\Delta(\varphi_1,\ldots,\varphi_m)$ is also simple.
\end{lemma}
\begin{proof}
The result is true for any component of $H=(H_1,\ldots,H_r)$ by
Lemma \ref{lema:simples}. Furthermore, if the statement is true
for some $\Delta$ then it is true for any divided difference
$\bar{\Delta}$ obtained by one of the formulas D1-D6, again by Lemma~\ref{lema:simples}.
\end{proof}
\begin{lemma}\labe{lema5}
Let $\Delta(\xx)\in\DD^m$ be a divided difference and denote by
$\widehat{\Delta}$ its Taylor series at $0\in\R^m$. Let
$\hat{\alpha}(x)=(\hat{\alpha}_1(x),\ldots,\hat{\alpha}_m(x))\in\R[[x]]^m$
be a vector of formal power series in a single variable with
$\hat{\alpha}(0)=0$. Then there are natural numbers $d_0,n_0\geq
0$ such that, for any $d>d_0$, there exists $L\in\N$, an analytic
function $C=C(x,\zz,\ww)$ in a neighborhood of $0\in\R^{1+m+rn_0}$
with $\partial\Delta/\partial\xx\,(0)=\partial
C/\partial\zz\,(0,0,0)$ and polynomials
$Q_1(x),\ldots,Q_{n_0}(x)\in\R[x]$ with $Q_j(0)=0$ such that
\begin{equation}\labe{eq:dd-formal}
T_d\left(\widehat{\Delta}(\hat{\alpha}(x))\right)=
C\left(x,T_d\hat{\alpha}\,(x),\{T_L\widehat{H}(Q_j(x))\}_{1\leq j\leq
n_0}\right).
\end{equation}
(Recall the notation $J_d\phi$ for the $d$-jet of $\phi$ at the
origin and the definition $T_d\phi(x)=(\phi(x)-J_d\phi(x))/x^d$).
\end{lemma}
The numbers $n_0$, $d_0$ (and also $C,\,Q_l,L$) do not depend
continuously (with respect to the topology of formal series) upon
$\hat \alpha(x)$, as can be seen in the proof. This means that
these objects can not be determined from a finite number of coefficients
of $\hat \alpha(x)$.

It is a tempting idea to try to prove the lemma using recursion, i.e.\ the
properties D1--D6; the problematic formula is D3. Any recursive proof
we found, however, is no simpler than the subsequent technical one.
\begin{remark}\labe{rem:dd-analitica}{\em
Before beginning the proof of Lemma \ref{lema5}, note that the result is
true (with $d_0=n_0=0$)
if we replace the divided difference $\Delta$ by an analytic
function $D\in\R\{\xx\}$: the analytic function
$C(x,\zz)=x^{-d}[D(J_d\hat\alpha(x)+x^d\zz)-J_d(D(\hat\alpha(x)))]$
verifies the statement. Moreover,
if the lemma holds for $\Delta_1$ and $\Delta_2$ then, chosing a suitable $d_0$,
it also holds for $\Delta_1+\Delta_2$. Therefore, given any
integer $M$, we can suppose that $J_M\Delta=0$.}
\end{remark}
\begin{proof}
We use formula (\ref{eq-dd}) of Lemma
\ref{lema3}. There might be a problem if $Q(\hat{\alpha}(x))$
vanishes identically. In any case, since the polynomial
$Q(\xx)$ is not identically zero, we can find some
$\tau=(\tau_1,\ldots,\tau_m)\in\N^m$ of minimal length
$|\tau|=\tau_{1}+\cdots+\tau_{n}$ such that
$\partial^{\tau}Q/\partial\xx^{\tau}(\hat{\alpha}(x))\not\equiv
0$. For convenience, put $P_0(\xx)=\partial^{\tau}Q/\partial\xx^{\tau}(\xx)$
and denote
$\hat{\beta_j}(x)=P_j(\hat{\alpha}(x))$ for $j=0,\ldots,n$. Up to
some permutation, we can suppose that there exists $0 \leq n_0\leq n$
such that $\hat{\beta}_j(x)\not\equiv 0$ if and only if
$j\in\{0,\ldots,n_0\}$. Let $\nu_j={\rm val}(\hat{\beta_j}(x))$
for $j\in\{0,\ldots,n_0\}$ and put $d_0={\rm
max}(\nu_0,(p+1)\nu_1,\ldots,(p+1)\nu_{n_0})$. Given $d>d_0$,
consider $M>d+\nu_0+|\tau|$. We can suppose, following
Remark \ref{rem:dd-analitica}, that $J_M\Delta=0$.

We consider the $\tau$-th derivative of equation (\ref{eq-dd})
\begin{equation}\labe{eq-derivadadd}
{\displaystyle \sum_{\tau'\leq\tau} }\left(
\begin{array}{c}
  \tau \\
  \tau' \\
\end{array}\right)
\frac{\partial^{\tau'}Q}{\partial\xx^{\tau'}}(\xx)
\frac{\partial^{\tau-\tau'}\Delta}{\partial\xx^{\tau-\tau'}}(\xx)=
\widetilde{B}\left(\xx,\{H^{(l)}(P_j(\xx))\}_{0\leq l\leq N+|\tau|,1\leq
j\leq n}\right)
\end{equation}
where $\widetilde{B}$ is some analytic function and by
$\tau'\leq\tau$ we mean that, if  $\tau'=(\tau'_1,\ldots,\tau'_m)$
then $0\leq\tau'_s\leq\tau_{s}$ for any $s$. Since
$x^{p+1}H'(x)=A(x,H(x))$, we deduce by induction
that $x^{(p+1)l}H^{(l)}(x)$ are analytic functions of $(x,H(x))$ for
all $l\in\N$. This fact, together with the
formula $H\,(x)=J_L H\,(x)+x^LT_{L}H\,(x)$ for any $L$,
leads to the property that there exists
$L\in\N$ and polynomials $R_{sj}(\xx)$ with ${\rm val}(R_{sj})\geq
M-|\tau|+1$ such that the right hand side of equation
(\ref{eq-derivadadd}) is of the form
\begin{equation}\labe{eq:btilde-bbarra}
\widetilde{B}\left( \xx,\{H^{(l)}(P_j(\xx))\}\right)= \overline{B}\left(
\xx,\{R_{sj}(\xx)\,T_{L}H_s(P_j(\xx))\}_{1\leq s\leq r,1\leq j\leq n}\right),
\end{equation}
with $\overline{B}$ analytic. Inserting $\xx=\hat{\alpha}(x)$ into equation
(\ref{eq-derivadadd}) we obtain from the minimality of $|\tau|$:
\begin{equation}\labe{eq:qtau-bbarra}
\hat{\beta}_0(x)
 \Delta(\hat{\alpha}(x))=
\overline{B}\left(\hat{\alpha}(x),\{R_{sj}(\hat{\alpha}(x))\,
T_{L}H_s(P_j(\hat{\alpha}(x)))\}_{1\leq s\leq r,1\leq j\leq
n_0},0\right).
\end{equation}
(Recall that for $j>n_0$ we have $P_j(\hat{\alpha}(x))\equiv 0$). Write for
$j=0,\ldots,n$
$$
P_j(J_{d}\hat{\alpha}(x)+x^{d}\zz)=J_d\hat{\beta}_{j}(x)+x^{d}D_{jd}(x,\zz)
$$
where $D_{jd}$ is some polynomial in the variables
$(x,\zz)\in\R^{1+m}$. By our choice of $d_0$ we can write for
$j=1,\ldots,n_0$:
\begin{equation}\labe{eq:p0-pj}
\begin{array}{l}
P_0(J_{d}\hat{\alpha}(x)+x^{d}\zz)=
J_d\hat{\beta}_{0}(x)\,(1+x^{d-\nu_0}\widetilde{D}_{0d}(x,\zz)),\\
P_j(J_{d}\hat{\alpha}(x)+x^{d}\zz)=
J_d\hat{\beta}_{j}(x)+(J_d\hat{\beta}_{j}(x))^{p+1}\widetilde{D}_{jd}(x,\zz),
\end{array}
\end{equation}
where $\widetilde{D}_{jd}$ are again analytic functions. Lemma
\ref{morfismo} (applied to $T_LH(x)=(H(x)-J_LH(x))/x^L$ which satisfies a
system of  ordinary differential equations analogous to (\ref{sistema}))
implies that for $s=1,\ldots,r$
and $j=1,\ldots,n_0$ there exists an analytic function $C_{sjdL}$
such that
$$
(T_LH_s)\,(P_j(J_d\hat{\alpha}(x)+x^d\zz))=C_{sjdL}\left(
x,\widetilde{D}_{jd}(x,\zz), (T_LH)(J_d\hat{\beta}_j(x))\right)
$$
On the other hand, if $j>n_0$ then $J_d\hat{\beta}_j(x)\equiv 0$ and
\begin{equation}\labe{eq:djd-nula}
D_{jd}(x,T_d\hat{\alpha}(x))= x^{-d}(P_j(J_d\hat{\alpha}(x)+x^dT_d\hat{\alpha}(x)))= x^{-d}P_j(\hat{\alpha}(x))\equiv 0. \end{equation}
We insert $\xx=J_d\hat{\alpha}(x)+x^d\zz$ into $\overline{B}$ and
obtain
\begin{eqnarray}\labe{eq:bbarra-bchec}
\overline{B}\left(
J_d\hat{\alpha}(x)+x^d\zz,\{R_{sj}(J_d\hat{\alpha}(x)+x^d\zz)\,
(T_LH_s)(P_j(J_d\hat{\alpha}(x)+x^d\zz))\}\right)=\\
\check{B} \left( x,\zz,\{x^{M_0}(T_LH)\,(J_d\hat{\beta}_j(x))\}_{j=1}^{n_0},\,
\{x^{M_0}(T_LH)\,(x^dD_{jd}(x,\zz))\}_{j=n_0+1}^{n}\right) \nonumber
\end{eqnarray}
with $M_0=M-|\tau|+1$ and some analytic function
$\check{B}$. Define for sufficiently small $\ww\in\R^{rn_0}$
$$
C(x,\zz,\ww):=\frac{\check{B} \left( x,\zz,\{x^{M_0}\ww\},0\right)}
{x^d\,J_d\hat{\beta}_{0}(x)(1+x^{d-\nu_0}\widetilde{D}_{0d}(x,\zz))}.
$$
Let us first show that $C$ is an analytic function. Since the
denominator in the quotient above has valuation $d+\nu_0<M_0-1$ it suffices
to show that the numerator is divisible by $x^{M_0-1}$.
We compute this numerator modulo $x^{M_0-1}$ from
(\ref{eq-derivadadd}), (\ref{eq:btilde-bbarra}) and
(\ref{eq:bbarra-bchec}):
$$
\check{B}(x,\zz,0,0) \equiv {\displaystyle \sum_{\tau'\leq\tau}} \left(
\begin{array}{c}
  \tau \\
  \tau' \\
\end{array}\right)
\frac{\partial^{\tau'}Q}{\partial\xx^{\tau'}}(J_d\hat{\alpha}(x)+x^d\zz)
\frac{\partial^{\tau-\tau'}\Delta}
{\partial\xx^{\tau-\tau'}}(J_d\hat{\alpha}(x)+x^d\zz).
$$
Since $J_M\Delta=0$, we obtain
$${\rm
val}_x\left(\frac{\partial^{\tau-\tau'}\Delta}{\partial\xx^{\tau-\tau'}}
(J_d\hat{\alpha}(x)+x^d\zz)\right)\geq M-(|\tau|-|\tau'|)\geq
M-|\tau|=M_0-1
$$
and thus $\check{B}(x,\zz,0,0) \equiv 0$ ${\rm mod}\, x^{M_0-1}$. Set
$Q_j(x)=J_d\hat{\beta}_j(x)$ for $j=1,\ldots,n_0$. We
use equations (\ref{eq:p0-pj}), (\ref{eq:djd-nula}),
(\ref{eq:bbarra-bchec}) and (\ref{eq:qtau-bbarra}) and that
$\hat{\alpha}(x)=J_d\hat{\alpha}(x)+x^dT_d\hat{\alpha}(x)$ and calculate:
$$
\begin{array}{l}
\!\!\!\!\!\!\!\!\!x^d\,\hat{\beta}_0(x)\,
C(x,T_d\hat{\alpha}(x),\{T_L\widehat{H}(Q_j(x))\}_j)=\\
=\check{B}(x,T_d\hat{\alpha}(x),\{x^{M_0}T_L\widehat{H}(J_d\hat{\beta}_j(x))\},0)=\\
=\overline{B}(\hat{\alpha}(x),\,\{R_{sj}(\hat{\alpha}(x))\,
(T_L\widehat{H}_s)(P_j(\hat{\alpha}(x)))\})=\\
=\hat{\beta}_0(x)\widehat{\Delta}(\hat{\alpha}(x)).
\end{array}
$$
This immediately yields (\ref{eq:dd-formal}).
Finally $\partial\Delta/\partial\xx\,(0)=\partial
C/\partial\zz\,(0,0,0)=0$ since $J_M\Delta=0$ and the function $C$
is divisible by $x$.
\end{proof}

\vspace{.3cm} {\em Proof of Proposition \ref{prop:a1=simples}\/}.-
Let $\beta\in\AAA^1$ and suppose that $\beta(0)=0$. By
Lemma \ref{lema2} there exists a system $\Phi(x,\yy)=0$ of $n$
implicit equations in $\AAA$ (here $x$ denotes a single variable and
$\yy\in\R^n$) with $\partial\Phi/\partial\yy\,(0)$ invertible,
such that each row of $\Phi$ is either analytic or of the form
$y_{l_i}-\Delta_i(x,\yy)$ for some $\Delta_i\in\DD^{1+n}$ and
such that its solution
$\yy=\phi=(\phi_1,\ldots,\phi_n)$ satisfies $\phi_n=\beta$. By Remark
\ref{rem:dd-analitica}, both kinds of rows satisfy the statement
of Lemma \ref{lema5} for the vector of formal series
$\hat{\alpha}(x)=(x,\hat{\phi}(x))\in\R[[x]]^{1+n}$, $\hat{\phi}$
being the Taylor series of $\phi$ at $x=0$. Taking the
maximum of the numbers $d_0$ for any of the rows and combining
the corresponding polynomials $Q_j$, we see that there exists an
analytic (vectorial) function $C(x,\zz,\ww)$ and $d\in\N$ such that
$\partial C/\partial\zz\,(0)=\partial\Phi/\partial\yy\,(0)$
is invertible and
\begin{equation}\labe{eq:Phi-formal}
C\left(
x,T_d\hat{\phi}(x),\{T_L\widehat{H}(Q_j(x))\}_j \right)=
T_d\left(\widehat{\Phi}(\hat{\alpha}(x))\right)\equiv0.
\end{equation}
The solution of the analytic equation $C(x,\zz,\ww)=0$ with
respect to $\zz$ is an analytic function $f(x,\ww)$. The germ
$\varphi(x)=f(x,\{T_LH(Q_j(x))\}_j)$ is a simple function and its
Taylor series $\hat{\varphi}(x)$ satisfies
$C(x,\hat{\varphi}(x),\{T_L\widehat{H}(Q_j(x))\}_j)\equiv 0$. By
the uniqueness of solutions of implicit equations,
$\hat{\varphi}(x)=T_d\hat{\phi}(x)$. Now, define the function
$$\widetilde{\varphi}(x)=J_d\phi(x)+x^d\varphi(x)=
(\widetilde{\varphi}_1(x),\ldots,\widetilde{\varphi}_n(x)).
$$ It is a vector of
simple functions whose Taylor series
$\widehat{\widetilde{\varphi}}(x)$ coincides with $\hat{\phi}(x)$.
Hence
$$
\widehat{\Phi}(x,\widehat{\widetilde{\varphi}}(x))=
\widehat{\Phi}(x,\hat{\phi}(x))\equiv 0.
$$
On the other hand, $\widehat{\Phi}(x,\widehat{\widetilde{\varphi}}(x))$ is the
Taylor series of $\Phi(x,\widetilde{\varphi}(x))$. This
last vector belongs to $(\SC_H)^n$ by Lemma \ref{lema4}. Due to
the quasi-analyticity of the class $\SC_H$, the germ of $x\mapsto
\Phi(x,\widetilde{\varphi}(x))$ is equal to $0$ and thus
$\phi=\widetilde{\varphi}\in(\SC_H)^n$. In particular
$\phi_n=\beta\in\SC_H$ is a simple function.\hfill{$\Box$}

\subsection{Theorem of the Complement for $\AAA$}

As a consequence of the quasi-analyticity results of the preceding subsection, the
theory developed in \cite{Rol-Spe-Wil} allows to complete the
proofs of Theorem \ref{main11} and Proposition \ref{H-conjuntos}.

Let $\AAA=\AAA_H$ be the smallest class of germs of $\CC^\infty$ functions
satisfying properties A1--A4. We define the $\AAA$-semianalytic and the
$\AAA$-subanalytic sets in the usual way. Given $a=(a_1,\ldots,a_m)\in\R^m$ and
$\delta>0$ denote $I_{a,\delta}=
[a_1-\delta,a_1+\delta]\times\cdots\times[a_m-\delta,a_m+\delta]$. A set
$A\subset I_{a,\delta}$ is called {\em $\AAA$-basic} if it is of the form
$$\{\xx\in I_{a,\delta}
\;/\;F(\xx)=0,G_{1}(\xx)>0,\ldots,G_{n}(\xx)>0\}
$$ for some $\AAA$-analytic
functions $F,G_{j}$ in a neighborhood of $I_{a,\delta}$. A set $A\subset\R^m$
is said to be {\em $\AAA$-semianalytic at the point $a\in\R^m$} if there exists
$\delta>0$ such that $A\cap I_{a,\delta}$ is a finite union of $\AAA$-basic
sets. We say that $A\subset\R^m$ is $\AAA$-semianalytic if it is
$\AAA$-semianalytic at any point of $\R^m$. In order to obtain finiteness of
the number of
connected components we restrict our considerations to {\em global $\AAA$-semianalytic
sets} $A\subset\R^m$, defined by the condition that $\rho_m(A)\subset I^m$ is a
$\AAA$-semianalytic set, where
$\rho_m\xunoxm=(x_1/\sqrt{1+x^{2}_{1}},\ldots,x_m/\sqrt{1+x^{2}_{m}})$.
Finally, a set $B\subset\R^m$ is said to be {\em (global) $\AAA$-subanalytic}
if there exists $n\geq m$ and a global $\AAA$-semianalytic set $A\subset\R^n$
such that $B=\pi(A)$, where $\pi:\R^n\rightarrow\R^m$ is the projection onto
the first $m$ coordinates. The following theorem is one of
the main results of the
article \cite{Rol-Spe-Wil}:
\begin{theorem}[Theorem of the Complement for $\AAA$] \labe{teorema-complementario}
Assume that $\AAA$ is a\linebreak quasi-analytic class satisfying A1--A4. Then, given any
$\AAA$-subanalytic set $A\subset\R^m$, the complement
$\R^m\setminus A$ is also $\AAA$-subanalytic and $A$ has a finite
number of connected components.
\end{theorem}
For completeness, let us briefly recall the idea of the proof
of Theorem \ref{teorema-complementario}. The quasi-analyticity of
the class $\AAA$ makes it possible to associate a finite {\em multiplicity}
to any non zero element $\phi\in\AAA^m$. Using the properties
A1--A4 of this class and a finite number of blowups and
ramifications, we can write $\phi$ in a {\em normal form}
$\phi(\xx)=m(\xx)u(\xx)$, where $m(\xx)$ is a monomial in the
coordinates $\xx$ and $u$ is a unit in $\AAA^m$. Normalizing all
the germs involved in the description of some $\AAA$-subanalytic set
$A$, one proves the finiteness of the number of connected
components of $A$, as well as the property that $A$ has a
dimension and that its boundary is included in some
$\AAA$-subanalytic set of dimension strictly smaller than that of
$A$. The property of the complement is a classical consequence of this property.%
\vspace{.5cm}

\noindent {\em Proof of Theorem \ref{main11}}. By Theorem
\ref{QA}, we can assume that the Theorem of the Complement for
$\AAA$ is true. Denote by $\SC\AAA$ the collection of global
$\AAA$-subanalytic sets. It is routine to check that $\SC\AAA$
is closed under finite unions and intersections, linear
projections and cartesian products. Also, it contains the graph of
any restricted analytic function as well as the graph of the
extension $\widetilde{H}$ of the solution $H$. It is closed under
taking complements by the first part of Theorem
\ref{teorema-complementario}. Thus, any definable set in the
structure $\R_{an,H}$ belongs to $\SC\AAA$ and $\R_{an,H}$ is
o-minimal by the second part of that theorem. On the other hand,
denote by $\mathcal{T}$ the smallest collection of sets satisfying
S1-S5 for $\FF=\widetilde{\FF}_{an}\cup\{\widetilde{H}\}$ without
the condition of taking complements and let us show that
$\SC\AAA\subset\mathcal{T}$. This will prove that the elements of
$\SC\AAA$ are exactly the definable sets in $\R_{an,H}$ and that
the structure $\R_{an,H}$ is model-complete. We only have to show
that the $\AAA$-basic sets (at the origin, for instance) are
contained in $\mathcal{T}$. Consider
the collection $\widetilde{\AAA}=(\widetilde{\AAA}^m)$ of
subalgebras of $\CC^\infty$ germs at $0\in\R^m$ having a
representative in some neighborhood of 0 whose graph is contained in
$\mathcal{T}$. We show that $\AAA\subset\widetilde{\AAA}$ by
verifying that $\widetilde{\AAA}$ satisfies properties A1--A4.
This proves that any $\AAA$-basic set of the type $\{\xx\in
I_{0,\delta}/F(\xx)=0\}$ belongs to $\mathcal{T}$. On the other
hand, any $\AAA$-basic set of the type $\{\xx\in
I_{0,\delta}/G(\xx)>0\}= \{\xx\in I_{0,\delta}/\exists
(y,z)\in\R\times\R^*\mbox{ with }G(\xx)=y \mbox{ and } y=z^2\}$ is
a projection of $({\rm graph}\,G\times\R)\cap X$ where $X$a is some
semialgebraic set. It also belongs to $\mathcal{T}$ and the proof is complete.

\vspace{.3cm}
 {\em Proof of Proposition \ref{H-conjuntos}}.- Denote by $\SC H$
the collection of $H$-subanalytic sets and let us prove that $\SC
H=\SC\AAA$. Clearly $\SC H\subset\SC\AAA$. On the other hand, it
is easy to check that $\SC H$ is closed under finite unions and
intersections, linear projections and cartesian products. We need
only to prove that the $\AAA$-basic sets belong to $\SC H$.
Consider the collection $\widetilde{\AAA}=(\widetilde{\AAA}^m)$ of
subalgebras of $\CC^\infty$ germs at $0\in\R^m$ having a
representative in some neighborhood of 0 whose graph is
$H$-subanalytic. The inclusion $\AAA\subset\widetilde{\AAA}$ follows by
verifying properties A1--A4 for $\widetilde{\AAA}$ and we conclude as
in the previous proof. \hfill{$\Box$}

\begin{remark}{\em
It is interesting to note that the collection ${\cal A}_H$ of algebras introduced below theorem \ref{main11} does {\em not} allow (in general) to
factor regular elements as in the Weierstrass preparation theorem which holds e.g.\ for the algebras of germs of functions analytic at the origin or
for those introduced by \cite{vdD-Spe} in the context of summability. Thus we have to rely on a result \cite{Rol-Spe-Wil} that uses desingularization
techniques not needing the Weierstrass preparation.

Precisely, consider any strictly increasing solution $H:[0,a]\ra\R,\ a>0$ of the equation $x^3 H'= (1-x^2)H-x$. As its formal solution $\hat
H(x)=\sum_{n\geq0}2^n\,n!\,x^{2n+1}$ is odd, $H$ can be extended to an odd solution $H:[-a,a]\ra\R$. Denote by ${\cal A}^m,m=1,2,...$ the algebras
generated by the germs of functions analytic at the origin and the germ of $H$ having the properties A1-A4 (see below theorem \ref{main11}). As the
formal solution $\hat H$ is divergent, the hypotheses of theorem~\deuxbis{} are satisfied. As discussed in example \ref{ex:curvas-pfaffianas}, if the
Stokes coefficients vanished then there would be no singular direction and the formal series would have to converge. Thus $\hat H$ is strongly
analytically transcendental (SAT).

Consider now the inverse function $I:[0,H(a)]\ra\R$ of $H$ which by proposition \ref{prop:a1=simples} is a simple function (this can also be shown
directly) and put $F(x,y)=I(y)^2+x$. This clearly is a function of ${\cal A}^2$ and it is regular as $F(0,y)=I(y)^2=y^2+{\cal O}(y^3)$, but {\em it
does not have a factor $y^2-G_1(x)y+G_2(x)$ that is polynomial in $y$ and belongs to ${\cal A}^2$} as we will show below.

First of all, the formal series $\hat F(x,y)=\hat I(y)^2+x$ can be factored $\hat F(x,y)=(y^2-\hat G_1(x)y+\hat G_2(x))\,\hat U(x,y)$ where $\hat
G_1,\hat G_2,\hat U$ are formal series and $\hat U(x,y)=1+...$ is a unit. The zeros of $\hat F$ are clearly $y=\hat H(\pm i\, x^{1/2})=\pm\hat
H(i\,x^{1/2})$ and thus $\hat G_1(x)=0,\ \hat G_2(x)=-\hat H(i\ x^{1/2})^2$. Putting $t=i\, x^{1/2}$ this yields $\hat G_2(-t^2)=-\hat H(t)^2$. This
implies that $\hat G_2(x)$ cannot be the Taylor series of a simple function $G_2=f(x,\{T_k H(P_l(x)),m=1,\ldots,n\})$ as in (\ref{eq-simples})
because an equation
$$f(-t^2,\{T_k \hat H(P_l(-t^2)),m=1,\ldots,n\})=-\hat H(t)^2$$
would contradict the (SAT); as $\hat H(x)$ is odd, $T_k \hat H(P_l(-t^2))$ can be reduced to some $\pm T_k \hat H(\tilde P_l(t^2))$, where $\tilde
P_l^{(\val \tilde P_l)}(0)>0$.

On the other hand, a factor $y^2-G_1(x)y+G_2(x)$ of $F(x,y)$ in ${\cal A}^2$ would have to have coefficients $G_1,G_2$ in ${\cal A}^1$, i.e.\ simple
functions according to proposition \ref{prop:a1=simples}. Therefore such a factor cannot exist.}
\end{remark}

%
\section{Quasi-analyticity of the algebra of simple functions}
\setcounter{theorem}{0}\setcounter{equation}{0}

This paragraph is devoted to the proof of Theorem \deuxbis{}. Let us
first show how this result implies Theorem \ref{main2}. This is an
immediate consequence of the following
\begin{lemma}\labe{lema:SATimplicaSQA}
Let $\widehat{H}(x)\in\R[[x]]^r$ be a formal solution of a system of
analytic ordinary differential equations of the form (\ref{sistema}). If $\widehat{H}(x)$ is
strongly analytically transcendental then any actual solution
$H:(0,\varepsilon]\rightarrow\R^r$ such that
$H(x)\sim\widehat{H}(x)$ as $\R^+\ni x\ra0$ is strongly quasi-analytic.
\end{lemma}
\begin{proof}
Let $f\in\R\{x,z_{11},\ldots,z_{rn}\}$ and
$P_1(x),\ldots,P_n(x)\in\R[x]$ be as in (SQA) and assume that
$f(x,\{T_k\widehat{H}_j\,(P_l(x))\}_{j,l})\equiv 0$. Let $\nu_l={\rm val} P_l>0$ and denote by
$Q_l(x)=J_{(p+1)\nu_l-1}\,P_l(x)$, a polynomial of valuation $\nu_l$
and degree $d_l<(p+1)\nu_l$.
Up to some permutation, we can assume that the set $\{Q_l|l=1,...,n\}=
\{Q_l|l=1,...,m\}$ with some $m\in\{1,...,n\}$ and {\em distinct}
$Q_l,\,l=1,...,m$. We can write
$$
P_l(x)=Q_l(x)+Q_l(x)^{p+1}u_l(x)
$$
where $u_l$ is some analytic function. By Lemma \ref{morfismo}
applied to $T_k H$ (which is a solution of a system of ordinary differential equations of the
form (\ref{sistema})), there exists another analytic function
$\widetilde{f}\in\R\{x,z_{11},\ldots,z_{rm}\}$ such that
$$f(x,\{T_k H_j\,(P_l(x))\}_{1\leq j\leq r,1\leq l\leq n})=
\widetilde{f}(x,\{T_k H_j\,(Q_l(x))\}_{1\leq j\leq r,1\leq l\leq m})\ .$$
Thus by assumption,
$\widetilde{f}(x,\{T_k\widehat{H}_j\,(Q_l(x))\}_{j,l})\equiv 0$.
By (SAT), we obtain $\widetilde{f}\equiv 0$ and we conclude that $f(x,\{T_k
H_j\,(P_l(x))\}_{j,l})\equiv 0$.
\end{proof}
\subsection{Background on elementary theory of summation and statement of the result}
Consider a system (\ref{sistema}) of analytic ordinary differential equations. We assume that
its linear part $A_0=\partial A/\partial\yy\,(0)$ has non zero
eigenvalues $\lambda_1,\ldots,\lambda_{r}$ satisfying the {\em
distinct argument condition}:
$$
\hspace{-3cm} (DA) \hspace{2cm} {\rm arg}(\lambda_i)\not\equiv{\rm
arg}(\lambda_j)\;\;{\rm mod}\;2\pi\Z\;\;\;\;\mbox{ if }i\neq j.
$$
As we have already mentioned in paragraph 1, condition (DA)
implies that there exists a unique formal power series solution
$\yy=\widehat{H}(x)=(\widehat{H}_1(x),\ldots,\widehat{H}_r(x))\in\R[[x]]^r$
with $\widehat H(0)=0$.
Its coefficients can easily be determined recursively due to the invertibility
of $A_0$\ .

\begin{remark}{\em
Moreover, although we do not use the result, there exists at least
one actual solution
$H=(H_1,\ldots,H_r):(0,\varepsilon]\rightarrow\R^r$ of
(\ref{sistema})
having the formal series $\widehat{H}(x)$ as its asymptotic
expansion as $\R^+\ni x\ra0$. When there are no eigenvalues on the imaginary
axis, the existence of $H$ is a consequence of the classical
theory of center manifolds (\cite{Hir-Pug-Shu,Car}). In the
general situation, we can proceed as in the fundamental theorem of
existence of complex valued solutions of (\ref{sistema}) \cite{Was}. See
also \cite{Bon-Dum} for the case $r=2$.}
\end{remark}

A general result about systems of ordinary differential equations
in the complex domain asserts that
formal solutions of systems (\ref{sistema}) are multisummable under very
weak conditions (see for instance \cite{Bra}); furthermore it is known that
the algebra of multisummable series is closed under composition. Instead of
using these facts as ``black box", we decided to present our
result using only some relatively elementary facts about summability of
divergent series. In our particular case where the linear part
$A_0$ is non singular we only need the following result (we use the
standard notation $S(\alpha,\beta;\rho)$, $\alpha,\beta\in\R$,
$\rho>0$ for the sector $\{z\in\C\,/\,| z|<\rho,\alpha<{\rm
arg}(z)<\beta\}$  in the complex plane of {\em opening}
$\beta-\alpha$):
\begin{lemma}\cite{Ram1}\labe{lema:H-sumable}
The vector $\widehat{H}(z)$ of formal series is $p$-summable. More
precisely, for any $\theta\in\R$, except for those which satisfy
$p\,\theta\equiv{\rm arg}\,(\lambda_j)\;{\rm mod}\,2\pi\Z$ for
some eigenvalue $\lambda_j$ of $A_0$, and for any open sector
$S_\theta=S(\theta-\frac{\pi}{2p}-\delta,\theta+\frac{\pi}{2p}+\delta;\rho)$
of opening slightly greater than $\pi/p$ (with $\delta,\rho$
sufficiently small) there exists a unique holomorphic function
$\widetilde{H}_\theta=
(\widetilde{H}_{\theta,1},\ldots,\widetilde{H}_{\theta,r}):S_\theta\rightarrow\C^r$
such that:
\begin{itemize}
\item[1)] $\widetilde{H}_\theta$ satisfies the (complexified) system of
equations (\ref{sistema}) and
\item[2)] For any $i$, $\widetilde{H}_{\theta,i}$ has the series
$\widehat{H}_i(z)$ as the asymptotic expansion on the sector
$S_{\theta}$. This asymptotic expansion is, furthermore, Gevrey of
order $\kappa=1/p$.
\end{itemize}
\noindent The function $\widetilde{H}_\theta$ is called the {\em
$p$-sum} of the series $\widehat{H}(z)$ along the ray
$d_\theta=\{z/{\rm arg}(z)=\theta\}$.
\end{lemma}
By definition (see \cite{Ram}), a bounded holomorphic function
$h:S\rightarrow\C$ on a sector $S=S(\alpha,\beta;\rho)$ has a {\em
Gevrey asymptotic expansion of order $\kappa$} (or a {\em
$\kappa$-Gevrey asymptotic expansion}) with right hand side
$\widehat{h}(z)=\sum_{k\geq 0} a_kz^k\in\C[[z]]$, and we write
$h(z)\sim_\kappa\widehat{h}(z)$, if for any $\eta>0$ there are
constants $K,A>0$ such that
$$
\left| h(z)-\sum_{k=0}^{N-1}a_kz^k\right|\leq K\,A^N\Gamma({N}{\kappa}+1)|
z|^N,
$$
for any $N\in\N$ and any $z\in S(\alpha+\eta,\beta-\eta;\rho-\eta)$.
Observe, that $h(z)\sim_\kappa 0+0z+...$ if and only if $h(z)$ is {\em exponentially
small of order $1/\kappa$}, i.e.\
for any $\eta>0$, there is a positive constant $a$ such that $|
h(z)\exp(a| z|^{-1/\kappa})|$ is bounded on
$S(\alpha+\eta,\beta-\eta;\rho-\eta)$. This follows easily from Stirling's Formula
by choosing $N$ as the integer closest to $| z| ^{-1/\kappa}$.

There are more algorithmical definitions of $p$-summability (see also e.g.\
\cite{Hor}). The function $\widetilde H_\theta(z)$ can be expressed as
{\em $p$-Borel-Laplace sum} of $\widehat H(z)$ in the direction $d_\theta$:
The {\em formal Borel transform} of order $p$, i.e.\ the series
${\bf H}(t)=\sum_{m\geq0}\frac1{\Gamma(1+m/p)}H_{m+1}t^m$ converges in
some neighborhood of $t=0$; it can be continued analytically to the
ray $\arg t=\theta$ (again named ${\bf H}(t)$) and has at most
exponential growth of order $p$; finally $\widetilde H_\theta$ is
the Laplace transform of order $p$ of ${\bf H}(t)$, i.e.\
$\displaystyle\widetilde H_\theta(z)=\int_{d_\theta}e^{-t^p/z^p}{\bf H}(t)
\frac{ p\,t^{p-1} }{z^{p-1}}\,dt\ \ .$

\vspace{.3cm} The excluded angles $\theta\in\{\frac{1}{p}({\rm
arg}\,(\lambda_j)+2\pi s)\;/\;j=1,\ldots,r,\,s\in\Z\}$ in the above
lemma are called {\em singular} and the corresponding
rays $d_\theta$ are called the {\em singular directions} (also
called the {\em anti-Stokes directions}) of $\widehat{H}(z)$. By
the hypothesis (DA), there are exactly $rp$ different singular
directions modulo $2\pi\Z$. Denote them by $d_{\theta_l}$, $l=0,\ldots,rp-1$ with
angles $\theta_l\in[0,2\pi)$. Up to reordering the eigenvalues
$\lambda_j$, we can choose the indices such that
$0\leq\theta_0<\cdots<\theta_{rp-1}$ and $p\,\theta_l\equiv {\rm
arg}(\lambda_j)$ ${\rm mod}\,2\pi\Z$ if and only if $l+1\equiv j$
${\rm mod}\,r$.

Two $p$-sums $\widetilde{H}_\theta$, $\widetilde{H}_{\theta'}$ of
$\widehat{H}(z)$ coincide in the intersection of their domains
if there is no
singular direction between $d_\theta$ and $d_{\theta'}$. We
obtain, by analytic continuation, holomorphic functions
$$
\widetilde{H}_l:\widetilde{S}_l=
S(\theta_l-\tfrac{\pi}{2p}+\delta,\theta_{l+1}+\tfrac{\pi}{2p}-\delta;\rho)
\rightarrow\C^r,\;\;\;l=0,\ldots,rp-1
$$
where $\delta,\rho>0$ are sufficiently small and where we put
$\theta_{rp}=\theta_0+2\pi$. Mostly, we will use some smaller sectors
$S_l=S(\theta_l-\varepsilon,\theta_{l+1}+\varepsilon;\rho)
\subset\widetilde{S}_l$ with some $\varepsilon>0$, so that the
family $\{S_l\}_l$ forms a {\em good covering} of the punctured disc
$\dot D_\rho\subset \C$,
i.e.\  $S_l\cap S_k\neq\emptyset$ if and only if
$| l-k|\leq 1$ ${\rm mod}\,rp$. In general, the functions
$\widetilde{H}_l$ cannot be continued analytically or change their
asymptotic behavior
beyond the rays of angles
$\theta_l-\frac{\pi}{2p},\theta_{l+1}+\frac{\pi}{2p}$,
$l=0,\ldots,rp-1$, called the {\em Stokes directions} of
$\widehat{H}(z)$. The so called {\em Stokes phenomenon} of the
series $\widehat{H}(z)$ is the description of the behavior of the
difference $\Delta_l(z)=\widetilde{H}_{l+1}(z)-\widetilde{H}_l(z)$
of two such consecutive functions. This difference is defined in the
intersection
$\widetilde{S}_{l,l+1}=\widetilde{S}_l\cap\widetilde{S}_{l+1}
=S(\theta_{l+1}-\frac{\pi}{2p}+\delta,\theta_{l+1}+\frac{\pi}{2p}-\delta;\rho)$,
a sector of opening slightly smaller than $\pi/p$. It satisfies a
system of linear differential equations (in the complex domain)
\begin{equation}\labe{eq:bl}
z^{p+1}\frac{d\yy}{dz}=B_l(z)\yy,
\end{equation}
where $B_l(z)$ is the matrix of holomorphic functions on
$\widetilde{S}_{l,l+1}$ defined by
$B_l(z)=\overline{B}(z,\widetilde{H}_l(z),\widetilde{H}_{l+1})$ with
$\overline{B}(z,\yy_1,\yy_2)$ analytic and satisfying
$$
A(z,\yy_2)-A(z,\yy_1)=\overline{B}(z,\yy_1,\yy_2)(\yy_2-\yy_1).
$$
We see that all matrices $B_{l}(z)$ for $l=0,\ldots,rp-1$ have the same
$\frac{1}{p}$-Gevrey asymptotic expansion on the sector $\widetilde{S}_{l,l+1}$
equal to $\widehat{B}(z)=\overline{B}(z,\widehat{H}(z),\widehat{H}(z))=\partial
A/\partial\yy\,(z,\widehat{H}(z))$. The initial term $\widehat{B}(0)=A_0$ has
distinct eigenvalues. In this situation, the classical theory of linear
ordinary differential equations
(see \cite{Was}) asserts that there exists a fundamental matrix solution of
(\ref{eq:bl}) of the form
\begin{equation}\labe{eq:Yl}
Y_l(z)=G_l(z)\,\exp(Q(z))\,z^{J},
\end{equation}
where:
\begin{enumerate}
\item $G_l(z)$ is a matrix of holomorphic functions on $\widetilde{S}_{l,l+1}$.
There exists a formal series $\widehat{G}(z)$ such that, for every $l$, the
matrix $G_l(z)$ has $\widehat{G}(z)$ as its $\frac{1}{p}$-Gevrey
asymptotic expansion in $\widetilde{S}_{l,l+1}$, moreover
${\rm det}(\widehat{G}(0))\neq 0$;
\item $J={\rm diag} (\alpha_1,\ldots,\alpha_r)$ is a constant diagonal matrix
and
\item  $Q(z)={\rm diag} (q_1(z),\ldots,q_r(z))$ is a diagonal matrix where
$q_j(z)=-\frac{\lambda_j}{p}z^{-p}+\cdots\in\C[z^{-1}]$ are polynomials in the
variable $z^{-1}$ of degree $p$ without constant term.
\end{enumerate}
The
anti-Stokes or singular directions (respectively the Stokes directions)
are precisely the rays where the initial term
of some of the polynomials $q_j(z)$ is real negative (respectively purely
imaginary).
Denote the columns of the matrix $G_l(z)$ by $G_{lj}(z)$,
$j=1,\ldots,r$.
The particular solution $\Delta_l$ of (\ref{eq:bl}) defined above (\ref{eq:bl})
can be written as
\begin{equation}\labe{eq:deltal}
\Delta_l(z)=\sum_{j} c_{lj}G_{lj}(z)e^{q_{j}(z)}z^{\alpha_j},
\end{equation}
where $c_l=(c_{l1},\ldots,c_{lr})\in\C^{r}$ is some constant vector.

\begin{lemma}\labe{lema:clmu-nonulo}
Given $l\in\{0,\ldots,rp-1\}$ and $\mu=\mu(l)\in\{1,\ldots,r\}$ defined by
$l+1\equiv \mu$ ${\rm mod}\,r$, we have that $c_{lj}=0$ for every $j\neq \mu$.
\end{lemma}
\begin{proof}
This is a classical result. It is due to the fact that for every $j\neq\mu$,
the function $\exp(q_j(z))$ is exponentially large on some ray in
$\widetilde{S}_{l,l+1}$ while $\Delta_l$ remains bounded on that sector.
\end{proof}

The coefficient $\gamma_l=c_{l\mu(l)}$ is called the {\em Stokes multiplier} of
the solution $\widehat{H}(z)$ associated to the singular direction
$d_{\theta_{l+1}}$.

Now we restate Theorem \deuxbis{} in the following more precise form:
\begin{theorem}\labe{main2-bis}
Consider a system of (real) analytic ordinary differential equations (\ref{sistema}) such that the linear part $A_0$ satisfies condition (DA).
Assume, furthermore, that the Stokes multipliers of its unique formal solution $\widehat{H}(z)\in\C[[z]]^r$ satisfy the following condition: for any
$\mu\in\{1,\ldots,r\}$ there exists some $l\in\{0,\ldots,pr-1\}$ with $l+1\equiv \mu$ ${\rm mod}\,r$ such that $\gamma_l\neq 0$. Then, the series
$\widehat{H}(z)$ is strongly analytically transcendental in the following sense: \vspace{.3cm}

\!\!\!\!  (SAT) \;\;\; If $k\geq 0$, $n\geq 0$, an
analytic function $f\in\C\{z,z_{11},\ldots,z_{rn}\}$ with $f(0)=0$
and distinct real polynomials $P_1(z),\ldots,P_n(z)\in\R[z]$
with $\deg P_l < (p+1) \val P_l $ and $P_l^{(\val P_l)}(0)>0$ are given, then one has
$$
f(z,\{T_k\widehat{H}_j\,(P_l(z))\}_{j,l})\equiv 0\; \Longrightarrow\; f\equiv
0.
$$
\end{theorem}

We will denote by (SD) (for {\em singularity of all directions ${\rm
arg}(\lambda_j)$}) the condition on the Stokes multipliers in the above
statement. In the proof of Theorem \ref{main2-bis} we will need the following
reformulations of two well known results of the theory of Gevrey series (see
\cite{Mal}).
\begin{lemma}[Ramis-Sibuya]\labe{lema:Ramis-Sibuya}
Let $\kappa>0$. Let $h_l:V_l\rightarrow\C$, $l=0,\ldots,N$ be a family of
bounded holomorphic functions in sectors $V_l=S(\alpha_l,\beta_l;\rho)$ such
that $\{V_0,\ldots,V_N\}$ is a good covering of the punctured disc
$\{z\in\C\,/\,0<| z|<\rho\}$ (i.e.\ the only two by two non-empty
intersections are the consecutive ones $V_{l,l+1}=V_l\cap V_{l+1}$, where
$V_{N+1}=V_0$). Then the differences $h_{l+1}-h_l$ are {\em exponentially small
of order $\kappa$} on $V_{l,l+1}$ for $l=0,\ldots,N$ if and only if there
exists a formal series $\widehat{h}(z)\in\C[[z]]$ such that
$h_l(z)\sim_{1/\kappa}\widehat{h}(z)$ on $V_l$ for any $l$.
\end{lemma}
\begin{lemma}[Relative Watson's Lemma]\labe{lema:Watson-relativo}
Let $\kappa_2>\kappa_1>0$ and $\beta>\alpha+\pi/\kappa_1$. Consider a good
covering of $V=S(\alpha,\beta;\rho)$ by sectors $V_l=S(\alpha_l,\beta_l;\rho)$,
$l=1,\ldots,N$. Suppose that we have holomorphic functions $h_l:V_l\rightarrow\C$
satisfying the following property: there are some constants $a,b>0$ such that
$| h_l(z)\exp(a\, |z|^{-\kappa_1})|$ and $|
(h_{l+1}(z)-h_l(z))\exp(b\,|z|^{-\kappa_2})|$ are bounded on $V_l$ and
$V_l\cap V_{l+1}$, respectively for $l=0,\ldots,N$. Then, $h_l$ is
exponentially small of order $\kappa_2$ on $V_l$ for any $l$.
\end{lemma}

\subsection{Proof of Theorem \ref{main2-bis}}
 First, taking into account that
$T_k\widehat{H}(z)$ is the unique formal solution of a system of equations
analogous to (\ref{sistema}),
we can suppose without loss of generality that $k=0$. Let
$f\in\C\{z,z_{11},\ldots,z_{rn}\}$ and $P_j(z)\in\R[z]$ be given satisfying the
hypothesis of (SAT). Consider the subset
$\Lambda\subset\{1,\ldots,r\}\times\{1,\ldots,n\}$ of indices
$(i,j)$ such that $f$ effectively depends on the variable
$z_{ij}$, i.e.\ its derivative with respect to this variable does not
vanish identically. We prove
Theorem \ref{main2-bis} by induction on the number of elements
of $\Lambda$.  Put
$\widehat{K}(z)=(z,\{\widehat{H}_i(P_j(z))\}_{(i,j)\in\Lambda})$ in order to
simplify the notation.

If $\Lambda$ is empty then the result is trivial.

Suppose that ${\rm card}\,\Lambda=c\geq 1$ and that the statement has been
shown if ${\rm card}\,\Lambda=c-1$ . In order to prove the result we assume
that $f(\widehat{K}(z))\equiv 0$ but $f\not\equiv 0$ and we deduce a
contradiction.

First we deduce that for every $(k,l)\in\Lambda$, there
exists some $s$ such that $\partial^s f/\partial
z_{kl}^s\,(\widehat{K}(z))\not\equiv 0$. In fact, if this property is not true,
then we have that the series in the two variables $z,z_{kl}$
$$
f(z,z_{kl},\{\widehat{H}_i(P_j(z))\}_{(i,j)\in\Lambda\setminus\{(k,l)\}})=
\sum_{m=0}^{\infty}\frac{1}{m!}\frac{\partial^m f}{\partial
z_{kl}^m}(\widehat{K}(z))\, (z_{kl}-\widehat{H}_k(P_l(z)))^m
$$
is identically zero. Thus, any coefficient of this series in $z_{kl}$ is
identically zero. Writing $f$ as a series
$f(z,\zz)=\sum_{m=0}^{\infty}f_m(z,\{z_{ij}\}_{(i,j)\neq (k,l)})z_{kl}^m$ we
conclude, by the hypothesis of induction on ${\rm card}(\Lambda)$, that
$f_m\equiv 0$ for all $m$ and hence $f\equiv 0$, contrary to our
assumption.

\labe{lema:parcialesnonulas}
For any fixed $(i,j)\in\Lambda$, we can suppose without loss of generality that
$\frac{\partial f}{\partial z_{ij}}(\widehat{K}(z))\not\equiv 0$.
The couple $(i,j)$ will be chosen below.
This can be justified as follows: Suppose for simplicity that $(i,j)=(1,1)$ and consider the smallest integer $s\geq 1$ such that $\partial^s f/\partial z_{11}^s\,(\widehat{K}(z))\not\equiv 0$.
Denote $g=\partial^{s-1}f/\partial z_{11}^{s-1}$. This function satisfies
$g(\widehat{K}(z))\equiv 0$ and $\partial g/\partial
z_{11}\,(\widehat{K}(z))\not\equiv 0$.
Renaming $g$ to $f$, we can assume the desired conditions for $f$.

Given $j\in\{1,\ldots,n\}$, the vector series $\widehat{H}(P_j(z))\in\C[[z]]^r$
is $\frac{1}{p\nu_j}$-summable, $\nu_j$ being the valuation of the
polynomial $P_j(z)$; this follows easily from Ramis-Sibuya's Lemma.
Moreover, $\varphi\in\R$ is a singular angle for
$\widehat{H}(P_j(z))$ if and only if $\theta=\nu_j\varphi$ is a singular angle
for $\widehat{H}(z)$. The $\frac{1}{p\nu_j}$-sum of $\widehat{H}(P_j(z))$ along
a non singular direction $d_\varphi=\{z/{\rm arg}(z)=\varphi\}$ is given by
$\widetilde{H}_{\nu_j\varphi}(P_j(z)):
\widetilde{V}_{j,\varphi}\rightarrow\C^r$, defined in some sector
$\widetilde{V}_{j,\varphi}$ of opening greater than $\pi/p\nu_j$ bisected by
$d_\varphi$ (recall the definition of $\widetilde{H}_{\theta}$ in Lemma
\ref{lema:H-sumable}). Let $\Gamma$ be the set of all singular directions
of the series $\widehat{H}(P_j(z)),\,j=1,...,n$
and denote its elements by $d_{\varphi_l}$, $l=0,\ldots,N$ with
$0\leq\varphi_0<\cdots<\varphi_N<2\pi$. If $d_\varphi\not\in\Gamma$ then the
function
$$
F_\varphi(z)=f(z,\{\widetilde{H}_{\nu_j\varphi,i}(P_j(z))\}_{(i,j)\in\Lambda})
$$
is defined in the sector
$\widetilde{V}_\varphi=\bigcap\widetilde{V}_{j,\varphi}$ of opening slightly
greater than $\pi/p\bar{\nu}$ where $\bar{\nu}={\rm max}(\nu_1,\ldots,\nu_n)$.
In a similar way as was explained before for the series $\widehat{H}(z)$, for
any given angle $\varphi$ in the interval $]\varphi_k,\varphi_{k+1}[$,
$k=0,\ldots,N$ (with $\varphi_{N+1}=\varphi_0+2\pi$), $F_\varphi$ can be continued
analytically to a holomorphic function
$\widetilde{F}_k:\widetilde{V}_k\rightarrow\C$ in a sector of the form
$\widetilde{V}_k=
S(\varphi_k-\frac{\pi}{2p\bar{\nu}}-\delta,\varphi_{k+1}+
\frac{\pi}{2p\bar{\nu}}+\delta;\rho)$
for some $\delta,\rho>0$ sufficiently small. We consider subsectors
$V_k=S(\varphi_k-\varepsilon,\varphi_{k+1}+\varepsilon;\rho)\subset\widetilde{V}_k$,
$k=0,\ldots,N$ which form a good covering of a punctured disc at $0\in\C$ and the
restrictions $F_k=\widetilde{F}_{k}|_{V_k}:V_k\rightarrow\C$. Let $\nu={\rm
min}(\nu_1,\ldots,\nu_n)$ and let us show the two following properties to
obtain the desired contradiction (from now on, any statement is valid, without mentioning
this explicitly every time, if $\varepsilon$ and $\rho$ are considered
sufficiently small in the definition of the sectors $V_k$):

I) Each $F_k$ is { exponentially small on $V_k$ of order $\kappa>p\nu$} on
$V_k$.

II) There is at least one $k_0\in\{0,\ldots,N\}$ such that the difference
$F_{k_0+1}-F_{k_0}$ is exponentially small of order exactly (not larger than)
$p\nu$ on the intersection $V_{k,k+1}=V_k\cap V_{k+1}$ (with $V_{N+1}=V_0$).

\vspace{.5cm} {\em Proof of I)}. Consider new variables
$\bar{\zz}=\{\bar{z}_{ij}\}_{(i,j)\in\Lambda}$ and analytic functions
$D_{ij}(z,\zz,\bar{\zz})$ satisfying
\begin{equation}\labe{eq:incrementos-finitos}
f(z,\bar{\zz})-f(z,\zz)=\sum_{(i,j)\in\Lambda}
D_{ij}(z,\zz,\bar{\zz})(\bar{z}_{ij}-z_{ij}).
\end{equation}
For any $(i,j)\in\Lambda$ and any $k$ denote by $h_{ij,k}:V_k\rightarrow\C$ the
restriction of the function $\widetilde{H}_{\nu_j\varphi,i}(P_j(z))$ with some
$\varphi\in(\varphi_k,\varphi_{k+1})$. Let $D_{ij,k}:V_{k,k+1}\rightarrow\C$ be
the function
obtained by replacing the variables $z_{uv}$ and $\bar{z}_{uv}$ in
$D_{ij}(z,\zz,\bar{\zz})$ by $h_{uv,k}(z)$ and $h_{uv,k+1}(z)$ respectively.
This yields for $k=0,\ldots,N$ and $z\in V_{k,k+1}$
\begin{equation}\labe{eq:Fl-Fl}
F_{k+1}(z)-F_k(z)=\sum_{(i,j)\in\Lambda} D_{ij,k}(z)(h_{ij,k+1}(z)-h_{ij,k}(z))\ .
\end{equation}
Notice that $D_{ij,k}$ has on $V_{k,k+1}$ a
Gevrey asymptotic expansion of order $1/p\nu$ with right hand side $\partial f/\partial
z_{ij}(\widehat{K}(z))$. The factor $h_{ij,k+1}-h_{ij,k}$ is exponentially small
of order greater or equal to $p\nu_j$ on $V_{k,k+1}$ by Ramis-Sibuya's Lemma,
since $h_{ij,k}$ and $h_{ij,k+1}$ are two $\frac{1}{p\nu_j}$-sums of
the same series $\widehat{H}_i(P_j(z))$. Thus, $F_{k+1}-F_k$ is exponentially
small of order greater or equal to $p\nu$ on $V_{k,k+1}$. Again by
Ramis-Sibuya's Lemma, the functions $F_k:V_k\rightarrow\C$ have a (common) Gevrey
asymptotic expansion of order $1/p\nu$, by construction equal to
$f(z,\{\widehat{H}_i(P_j(z))\})\equiv 0$. Hence $F_k$ is exponentially small of
order greater or equal to $p\nu$ on $V_k$ .

Let us show that this order is
strictly greater than $p\nu$. First, note that $F_k$ is defined as $F_\varphi$
for any $\varphi$ in the interval $(\varphi_k,\varphi_{k+1})$ and that
$F_\varphi$ is defined at least in the sector $\widetilde{V}_\varphi$ bisected
by the ray $d_\varphi$ and whose opening does not depend on $\varphi$. Hence,
it suffices to prove
\begin{lemma}\labe{lema:ftheta}
If $d_\varphi\not\in\Gamma$ then $F_\varphi$ is exponentially small of order
(strictly) greater than $p\nu$ in $\widetilde{V}_\varphi$.
\end{lemma}
\begin{proof}
Fix $d_\varphi\not\in\Gamma$ and consider the set
of sectors $\{V_{k_1},\ldots,V_{k_2}\}$
that have a non empty intersection with the closed sector $\{z\,\mid\,|
z|<\rho,|\arg(z)-\varphi|\leq\frac{\pi}{2p\nu}\}$. The union
$V=V_{k_1}\cup\cdots\cup V_{k_2}$ is a sector of opening strictly greater than
$\pi/p\nu$ and contains $d_\varphi$. We can suppose that
$\widetilde{V}_\varphi\subset V$. Moreover, we can suppose that for any
$j\in\{1,\ldots,n\}$ such that $\nu_j=\nu$, the holomorphic function
$\widetilde{H}_{\nu\varphi}(P_j(z))$ is defined in the whole sector $V$.
Consider for any $k\in\{k_1,\ldots,k_2\}$ the holomorphic function
$\overline{F}_k:V_k\rightarrow\C$ defined by
$\overline{F}_k(z)=f(z,\{\overline{h}_{ij,k}(z)\}_{(i,j)\in\Lambda})$ where
$\overline{h}_{ij,k}$ is the suitable modification of the function $h_{ij,k}$
given by
$$
\overline{h}_{ij,k}(z)=\left\{%
\begin{array}{ll}
    \widetilde{H}_{\nu\varphi,i}(P_j(z)), & \hbox{if $\nu_j=\nu$;} \\
    h_{ij,k}(z), & \hbox{if $\nu_j>\nu$.} \\
\end{array}%
\right.
$$
Notice that the differences $\overline{h}_{ij,k+1}-\overline{h}_{ij,k}$ are
either zero or equal to $h_{ij,k+1}-h_{ij,k}$ when $\nu_j>\nu$, so they are
exponentially small of order strictly greater than $p\nu$ on $V_{k,k+1}$.
Using an equation analogous to  (\ref{eq:Fl-Fl}), the same is true for
$\overline{F}_{k+1}-\overline{F}_k$. By the Relative Watson's Lemma, the
functions $\overline{F}_k:V_k\rightarrow\C$ are also exponentially small of
order greater than $p\nu$. Finally, we compare $\overline{F}_k$ to
$F_\varphi$  using an equation similar to (\ref{eq:Fl-Fl}): there are
bounded functions $\overline{D}_{ij,k}(z)$ in $\widetilde{V}_\varphi\cap V_k$
with
$$
\overline{F}_k(z)-F_\varphi(z)=\sum_{(i,j)\in\Lambda}
\overline{D}_{ij,k}(z)(\overline{h}_{ij,k}(z)-\widetilde{H}_{\nu_j\varphi,i}(P_j(z))).
$$
Hence, $\overline{F}_k-F_\varphi$ is exponentially small of order strictly
greater than $p\nu$ on $\widetilde{V}_\varphi\cap V_k$ for any
$k\in\{k_1,\ldots,k_2\}$. The result follows.
\end{proof}

\vspace{.5cm} {\em Proof of II.-} We want to estimate the exponential growth of
the difference $F_{k+1}(z)-F_k(z)$, for some index $k$ chosen below. We use
equation (\ref{eq:Fl-Fl}). We first look at the second factor of each summand
in that expression. We have for any $(i,j)\in\Lambda$ and $z\in V_{k,k+1}$:
\begin{equation}\labe{eq:hnuj-hnuj}
h_{ij,k+1}(z)-h_{ij,k}(z)=\widetilde{H}_{\nu_j\overline{\varphi},i}(P_j(z))-
\widetilde{H}_{\nu_j\varphi,i}(P_j(z)),
\end{equation}
where $\varphi<\varphi_{k+1}<\overline{\varphi}$ are arbitrarily chosen but
sufficiently close to $\varphi_{k+1}$. Notice that if $\nu_j\varphi_{k+1}$ is
not one of the singular angles $\{\theta_0,\ldots,\theta_{rp-1}\}$ of the
series $\widehat{H}(z)$ then $\widetilde{H}_{\nu_j\varphi}$ and
$\widetilde{H}_{\nu_j\overline{\varphi}}$ coincide and hence
$h_{ij,k+1}-h_{ij,k}\equiv 0$.
 Otherwise, if $\nu_j\varphi_{k+1}=\theta_{l(k,j)+1}$ for some
 $l(k,j)\in\{0,\ldots,rp-1\}$ (recall that we set $\theta_{rp}=\theta_0+2\pi$) then we have
$h_{ij,k+1}(z)-h_{ij,k}(z)=\Delta_{l(k,j),i}(P_j(z))$ (the function
$\Delta_l=(\Delta_{l,1},\ldots,\Delta_{l,r})=\widetilde{H}_{l+1}-\widetilde{H}_l$
being defined in the preceding paragraph). Denote by $J_k$ the set of indices $j$
such that there exists $l(k,j)$ with
$\nu_j\varphi_{k+1}=\theta_{l(k,j)+1}$. Using equation
(\ref{eq:deltal}), Lemma \ref{lema:clmu-nonulo} and the definition
of the Stokes multipliers $\gamma_l$ we obtain for any $k$
\begin{align}\labe{eq:Fk-Fk-bis}
F_{k+1}&(z)-F_k(z)=  \\
&\sum_{(i,j)\in\Lambda,\,j\in J_k}
D_{ij,k}(z)\gamma_{l(k,j)}G_{l(k,j)\mu(k,j)}^{i}(P_j(z))
\exp(q_{\mu(k,j)}(P_j(z)))P_j(z)^{\alpha_{\mu(k,j)}}.\nonumber
\end{align}
Here $\mu(k,j)\in\{1,\ldots,r\}$ is defined by $\mu(k,j)\equiv l(k,j)+1$ ${\rm mod}\,r$ and $G_{l\mu}^{i}$ is the $i$-th entry of the column vector $G_{l\mu}$
of the matrix $G_l$ in (\ref{eq:Yl}).

We compute the real part of $q_{\mu(k,j)}(P_j(z))$:
\begin{equation}\labe{eq:partereal-exponenciales}
 \Re(q_{\mu(k,j)}(P_j(z)))=-\frac{|\lambda_{\mu(k,j)}|}{pa_{j}^{p}|
 z|^{p\nu_j}}\cos({\rm arg}(\lambda_{\mu(k,j)})-p\nu_j{\rm
 arg}(z))\,(1+O(| z|)),
\end{equation}
where here and in the sequel $a_j:=P_j^{(\nu_j)}(0)/\nu_j!>0$ with $\nu_j=\val P_j$
denotes the first non-zero coefficient of $P_j$.
Notice that if $z$ is in a small sector bisected by
$d_{\varphi_{k+1}}$ then $p\nu_j{\rm arg}(z)$ is close to ${\rm
arg}(\lambda_{\mu(k,j)})$ ${\rm mod}\,2\pi\Z$. Hence, the
exponential term in (\ref{eq:Fk-Fk-bis}) satisfies for $z\in
V_{k,k+1}$ and $\delta>0$ sufficiently small:
\begin{equation}\labe{eq:modulo-exponenciales}
|\exp(q_{\mu(k,l)}(P_j(z)))|=\exp(\Re(q_{\mu(k,l)}(P_j(z))))=O\left(\exp
\left(-\tfrac{|\lambda_{\mu(k,j)}|-\delta}{pa_{j}^{p}|
 z|^{p\nu_j}}\right)\right).
\end{equation}
\vspace{.2cm} Now we describe how to choose a suitable index $k_0$. Suppose for
instance that $\nu=\nu_1$. Up to a permutation of the indices, we can suppose
that $P_1,\ldots,P_{n_1}$ are the polynomials having valuation
$\nu$ and minimal first nonzero coefficient $A:=a_{1}=\cdots=a_{n_1}$
(i.e.\ if $j>n_1$ then either $\nu_j>\nu$ or $\nu_j=\nu$ and $a_{j}>A$).
Consider the subset
$I\subset\{0,\ldots,N\}$ of all indices $k$ for which
 the (unique)
singular direction $d_{\varphi_{k+1}}\in\Gamma$ in $V_{k,k+1}$ is
a singular direction of the series $\widehat{H}(P_1(z))$ (and
hence of any $\widehat{H}(P_j(z))$ for $j=1,\ldots,n_1$). With the
above notations, if $k\in I$ then $\{1,\ldots,n_1\}\subset
J_k$ and in this case $l(k):=l(k,1)=\cdots=l(k,n_1)$ is such that
$\theta_{l(k)+1}=\nu\varphi_{k+1}$. Define also
$\mu(k)\in\{1,\ldots,r\}$ by $\mu(k)\equiv l(k)+1$ ${\rm mod}\,r$
for any $k\in I$. Then equation (\ref{eq:Fk-Fk-bis}) yields with
(\ref{eq:modulo-exponenciales}) for $k\in I$
\begin{align}\labe{eq:Fk-Fk-bisbis}
F_{k+1}(z)-F_k(z)=
\gamma_{l(k)}&\sum_{j=1}^{n_1}E_{k,j}(z)\exp(q_{\mu(k)}(P_j(z)))
P_j(z)^{\alpha_{\mu(k)}}+\nonumber\\
& +O(\exp (-B| z|^{-p\nu})),
\end{align}
where $B$ is a certain constant
$>\frac{|\lambda_{\mu(k)}|}{pA^p}$ and, for simplicity, we have denoted
$$
E_{k,j}(z)=\sum_{i\,/\,(i,j)\in\Lambda}D_{ij,k}(z)G_{l(k,j)\mu(k,j)}^{i}(P_j(z)).
$$

\begin{lemma}\labe{lema:mu0}
We can assume, without loss of generality,
that there exists $\mu_0\in\{1,\ldots,r\}$ with the property
that $E_{k,1}(z)$ has a non trivial asymptotic expansion on $V_{k,k+1}$ for any
$k\in I$ with $\mu(k)=\mu_0$.
\end{lemma}
\begin{proof}
The asymptotic expansion of $E_{k,1}$ for $k\in I$ is given by
$$
\widehat{E}_{k,1}(z)=(\widehat{D}_{11,k}(z),\ldots,
\widehat{D}_{r1,k}(z))\,\widehat{G}(P_j(z))\,
(0,\ldots,\stackrel{(\mu(k))}{1},\ldots,0)^T
$$
where $\widehat{D}_{ij,k}(z)=\partial f/\partial
z_{ij}(z,\widehat{K}(z))\in\C[[z]]$ is the asymptotic expansion of $D_{ij,k}$.
As discussed in the beginning of this section, we can assume that the series $\partial
f/\partial z_{11}(z,\widehat{K}(z))\in\C[[z]]$ is not identically zero,
hence the vector $(\widehat{D}_{11,k}(z),\ldots,$ $
\widehat{D}_{r1,k}(z))$ is not trivial. We
conclude using the fact that $\widehat{G}(z)$ is a non singular matrix and the
surjectivity of the function $\mu:I\rightarrow\{1,\ldots,r\}$.
\end{proof}

Now we use the hypothesis (SD) of Theorem \ref{main2-bis}: choose $l_0\in\{0,\ldots,rp-1\}$ with $l_0+1\equiv\mu_0$ ${\rm
mod}\,r$ such that $\gamma_{l_0}\neq 0$ and consider
$k_0\in\{0,\ldots,N\}$ such that $\theta_{l_0+1}=\nu\varphi_{k_0+1}$. Then $k_0\in
I$ and $l_0=l(k_0)$, $\mu_0=\mu(k_0)$. We are now in the position to show that
$F_{k_0+1}-F_{k_0}$ is exponentially small of the exact order $p\nu$.

We use (\ref{eq:Fk-Fk-bisbis}) for  $k=k_0$. By Lemma
\ref{lema:mu0}, at least the coefficient $E_{k_0,1}$ has a non zero
asymptotic expansion on $V_{k_0,k_0+1}$. We can suppose, without loss of generality,
that this is
true for any coefficient $E_{k_0,j}$ for $j=1,\ldots,n_1$. In fact, if some
$E_{k_0,j}$ has a vanishing  asymptotic expansion then $E_{k_0,j}$ is
exponentially small of order $p\nu$ on $V_{k_0,k_0+1}$ (since this expansion is
Gevrey of order $1/p\nu$). This property, together with
(\ref{eq:modulo-exponenciales}), permits to include the corresponding summand
$E_{k_0,j}(z)\exp(q_{\mu_0}(P_j(z)))P_j(z)^{\alpha_{\mu_0}}$ in the remainder
term $O(\exp(-B| z|^{-p\nu}))$ in (\ref{eq:Fk-Fk-bisbis})
(maybe the constant $B$ has to be reduced but remains $>\frac{|\lambda_{\mu(k)}|}{pA^p}$).

We remark that we have not yet used the hypothesis
of Theorem \ref{main2-bis} that the degree $d_j$ of the
polynomial $P_j(x)$ satisfies $d_j<(p+1)\nu_j$. This hypothesis
(together with the fact that the polynomial
$q_{\mu_0}(z)\in\R[z^{-1}]$ is of degree $p$ and without constant
term) is essential in order to show the following property (its
proof is straightforward and left to the reader):
\begin{quote}
The principal parts at $z=0$ of the meromorphic series $q_{\mu_0}(P_j(z))$ for
$j=1,\ldots,n_1$ are two by two different.
\end{quote}
Thus, for any pair of distinct indices $j_1,j_2\in\{1,\ldots,n_1\}$, there
exists a complex number $\beta_{j_1j_2}\neq 0$ and an integer $s_{j_1j_2}>0$
such that
$$
q_{\mu_0}(P_{j_1}(z))-q_{\mu_0}(P_{j_2}(z))=
\beta_{j_1j_2}z^{-s_{j_1j_2}}(1+O(| z|)).
$$
We deduce that for almost all $\varphi\in\R$, except for those satisfying ${\rm
arg}(\beta_{j_1j_2})-s_{j_1j_2}\varphi-\frac\pi2\in\pi\Z$, the real part of
$q_{\mu_0}(P_{j_1}(z))-q_{\mu_0}(P_{j_2}(z))$ has a constant sign (positive or
negative) along $d_\varphi$ ( for $| z|$ sufficiently small depending on
$\varphi$). Up to a permutation of the indices, we can suppose that
$\Re(q_{\mu_0}(P_{1}(z)))>\Re(q_{\mu_0}(P_{j}(z)))$ for $j>1$ along infinitely
many rays $d_\varphi$ where $\varphi$ is arbitrarily close to
$\varphi_{k_0+1}$. Moreover, given such an angle $\varphi$, there exists
$b=b(\varphi)>0$ and $s>0$ (independent of $\varphi$) such that
$$
Re(q_{\mu_0}(P_{1}(| z| e^{i\varphi})))-\Re(q_{\mu_0}(P_{j}(| z|
e^{i\varphi})))>b| z|^{-s},
$$
for $| z|$ small. Due to the fact that $E_{k_0,j}$ has a non
trivial asymptotic expansion, this implies that the first term in
the summation in (\ref{eq:Fk-Fk-bisbis}) dominates
the remaining ones; i.e.
$$
\lim_{z\rightarrow 0,z\in d_\varphi}\frac{|
E_{k_0,j}(z)\exp(q_{\mu_0}(P_j(z)))P_j(z)^{\alpha_{\mu_0}}|} {|
E_{k_0,1}(z)\exp(q_{\mu_0}(P_1(z)))P_1(z)^{\alpha_{\mu_0}}|}=0,
$$
for $j=2,\ldots,n_1$. It also dominates the term $O(-B|
z|^{-p\nu})$ since $B>\frac{\mid\lambda_{\mu(k)}\mid}{pA^p}$. Thus we have
$$
F_{k_0+1}(z)-F_{k_0}(z)=
\gamma_{l_0}E_{k_0,1}(z)\exp(q_{\mu_0}(P_1(z)))P_1(z)^{\alpha_{\mu_0}}(1+O(z))
$$
as $z\rightarrow0$ restricted to $z\in d_{\varphi}$ for some $\varphi$ arbitrarily close to
$\varphi_{k_0+1}$. Finally, by (\ref{eq:partereal-exponenciales}), the function
$\gamma_{l_0}E_{k_0,1}(z)\exp(q_{\mu_0}(P_1(z)))P_1(z)^{\alpha_{\mu_0}}$ is
exponentially small of exact order $p\nu$ and thus the same is true for
$F_{k_0+1}-F_{k_0}$.

 \end{document}